\RequirePackage{ifpdf}
\ifpdf 
\documentclass[pdftex]{sigma}
\else
\documentclass{sigma}
\fi

\numberwithin{equation}{section}

\newtheorem{Theorem}{Theorem}[section]
\newtheorem{Corollary}[Theorem]{Corollary}
{\theoremstyle{definition}
\newtheorem{Example}[Theorem]{Example}
\newtheorem{Remark}[Theorem]{Remark}
}

\def\b1{\text{\bf 1}}
\def\B1{{\Bbb 1}}

\def\BA{{\Bbb A}}
\def\BC{{\Bbb C}}
\def\BF{{\Bbb F}}
\def\BH{{\Bbb H}}
\def\BN{{\Bbb N}}
\def\BP{{\Bbb P}}

\def\BZ{{\Bbb Z}}

\def\CC{{\cal C}}

\def\CL{{\cal L}}

\def\CO  {{\cal O}}

\def\CT{{\cal T}}

\def\dpar{\partial}

\def\fR{{\frak R}}

\def\hra{\hookrightarrow}

\def\isom{\buildrel\sim\over{=}}

\def\lra{\longrightarrow}

\begin{document}

\allowdisplaybreaks

\renewcommand{\thefootnote}{$\star$}

\renewcommand{\PaperNumber}{034}

\FirstPageHeading

\ShortArticleName{Homological Algebra and Divergent Series}

\ArticleName{Homological Algebra and Divergent Series\footnote{This paper is a
contribution to the Special Issue on Kac--Moody Algebras and Applications. The
full collection is available at
\href{http://www.emis.de/journals/SIGMA/Kac-Moody_algebras.html}{http://www.emis.de/journals/SIGMA/Kac-Moody{\_}algebras.html}}}

\Author{Vassily GORBOUNOV~$^\dag$ and Vadim SCHECHTMAN~$^\ddag$}

\AuthorNameForHeading{V.~Gorbounov and V.~Schechtman}

\Address{$^\dag$~Department of Mathematical Sciences, King's
College, University of Aberdeen,\\
\hphantom{$^\dag$}~Aberdeen, AB24 3UE, UK}
\EmailD{\href{mailto:vgorb@maths.abdn.ac.uk}{vgorb@maths.abdn.ac.uk}}

\Address{$^\ddag$~Laboratoire de Math\'ematiques Emile Picard,
Unversit\'e Paul Sabatier, Toulouse, France}
\EmailD{\href{mailto:schechtman@math.ups-tlse.fr}{schechtman@math.ups-tlse.fr}}

\ArticleDates{Received October 01, 2008, in f\/inal form March 04,
2009; Published online March 24, 2009}

\Abstract{We study some features of inf\/inite resolutions of Koszul
algebras motivated by the developments in the string theory
initiated by Berkovits.}

\Keywords{Koszul resolution; Koszul duality, divergent series}

\Classification{13D02; 14N99}



\medskip

\rightline{\it To Friedrich Hirzebruch on his 80-th anniversary, with admiration}

\renewcommand{\thefootnote}{\arabic{footnote}}
\setcounter{footnote}{0}

\section{Introduction}

\paragraph{1.1.} This article consists of two parts. The f\/irst part
is a simple exercise on Mellin transform. The second one is a review
on some numerical aspects of Koszul duality. An object which lies
behind the two parts is a Tate resolution of a commutative ring over
a f\/ield of characteristic zero.

Let us give some more details on the contents. In the f\/irst part we
develop the elegant ideas, due to physicists~\cite{BN}, which allow
to def\/ine the numerical invariants of projective varieties, doing a
regularization of some divergent series connected with their
homogeneous rings.

Let $R_0 = k[x_0,\ldots,x_N]$ be a polynomial algebra over a f\/ield
$k$,  $f_1, \ldots, f_p \in R_0$ homogeneous elements of degrees
$d_i = \deg f_i >0$ which generate the ideal $I = (f_1,\ldots,f_p)$.
Consider the quotient algebra $A = R_0/I$; it is graded
\[
A =  \mathop{\oplus}\limits_{j=0}^\infty  A^j.
\]
Geometrically the projection
$R_0 \lra A$ corresponds to a closed embedding
$i: X := \text{Proj}\, A \hra \BP^N := \text{Proj}\, R_0$; the projective
variety  $X$ is def\/ined in $\BP^N$ by the equations $f_1 = 0, \ldots,
f_p = 0$. The algebra $A$ is called the homogeneous ring of $X$ (it depends
on the embedding into the projective space).

Let us call {\it a semi-free resolution} of $A$ the following data:

$(i)$ An associative unital bi-graded $k$-algebra $R =
\mathop{\oplus}\limits_{i,j=0}^\infty  R_i^j$ (so $R_i^j\cdot R_l^m \subset
R_{i+l}^{j+m}$); the indexes~$i$ and~$j$ will be called the
homological degree and the polynomial degree respectively. We set
$R_i := \oplus_j  R_i^j$, $R^j := \oplus_i  R_i^j$. We introduce a
structure of a superalgebra on $R$ by def\/ining the parity to be
equal to the parity of the homological degree.

The multiplication has to be super-commutative, i.e.\ for $x\in R_i$,
$y\in R_l$ we have $xy = (-1)^{il}yx$.

$R$ should be equipped with a dif\/ferential $d:  R\lra R$, $d^2 = 0$, such that
$d(R_i^j) \subset R_{i-1}^j$ and for $x\in R_i$, $y\in R$,
\[
d(xy) = dx\cdot y + (-1)^ix\cdot dy.
\]

$(ii)$ $R_0$ should be equipped with a morphism of algebras
$\epsilon:  R_0\lra A$ such that  $\epsilon(R_0^j)\subset A^j$.

$(iii)$ For each $j\geq 0$, $R^j$ should be a resolution of $A^j$, i.e.\ the
complexes
\[
\cdots \lra R_2^j \buildrel{d}\over\lra R_1^j
\buildrel{d}\over\lra R_0^j \buildrel\epsilon\over\lra A^j\lra 0
\]
should be exact.

$(iv)$ If one forgets the dif\/ferential
$d$, $R$ should be a polynomial (super)algebra in homogeneous generators.
We suppose that in each polynomial degree one has a f\/inite number
of generators.

Semi-free resolutions exist for each $A$ but they are not unique
(sometimes they are called {\it Tate resolutions}).

\begin{Example} Let us consider the Koszul complex
\[
K(R_0;\text{\bf f}):=
k[x_0,\ldots,x_N]\otimes \Lambda\langle \xi_1\cdots \xi_M\rangle,
\]
where $\Lambda$ denotes the exterior algebra; the dif\/ferential is
def\/ined by the formula $d\xi_j = f_j$. The requirements $(i)$, $(ii)$,
$(iv)$ are fulf\/illed (the homological degree of $x_i$ is $1$). The
requirement $(iii)$ holds if\/f the sequence $\text{\bf f} =
\{f_1,\ldots, f_M\}$ is regular.

If this is the case, the Krull dimension of $A$ is given by a simple formula
\begin{gather}
\dim A = N+1 - M, \label{(1.1.1)}
\end{gather}
where $N+1$ is the number of $x_i$ and $M$ is the number of $\xi_j$;
one can say that $N+1 - M$ is the ``superdimension'' of the polynomial
ring $k[x_0,\ldots,x_N]\otimes \Lambda\langle \xi_1\cdots \xi_M\rangle$ $($the
number of even generators minus the number of odd generators$)$.
\end{Example}

In general one can consider a semi-free resolution of $A$ as a
natural replacement of the Koszul complex. One of the main goals of
this paper is {\it to propose an analog of the equation \eqref{(1.1.1)}
valid for not necessarily complete intersections}.

If our algebra $A$ is not a complete intersection then a semi-free
resolution will be inf\/inite. However, one can associate with such a
resolution a sequence of integer numbers and to show that  $\dim A$
can be written as a ``sum'' of these numbers; this will be the
required generalization of \eqref{(1.1.1)}. One puts the sum inside of the
quotation marks since we have here a divergent series, the sum of
which is calculated be means of a regularization, following  a
classical procedure of Riemann.

In fact, our formula is nothing else but the Mellin transform of the
classical theorem by Hilbert which says that $\dim A$ is equal to
the order of pole of the Hilbert series of $A$ at $1$. Moreover, we
give a similar interpretation for the other numerical invariants of
$A$ such as the degree. For details, see Sections \ref{section2.1}--\ref{section2.4}.

Except for the case of the complete intersection which is rather
seldom, there exists a wider class of algebras which admit a
remarkable explicit semi-free resolution. Namely, if $A$ is a~{\it
Koszul algebra} then one can take for $R$ the Chevalley cochain
complex of the graded Lie algebra~$L$ Koszul dual to $A$; one can
call it {\it the Koszul--Chevalley resolution}, cf.\ Section~\ref{section3.4} for
details.

As has been noted by many the notion of the Koszul duality is of
fundamental importance in Physics. For us a motivating reference is
the work of M.~Movshev and A.~Schwarz. It turns out it is closely
related to the ``gauge f\/ields-strings duality'' dear to Polyakov,
cf.~\cite{Pol}.

In the second part of this article we discuss a (probably) simplest
nontrivial case of a Koszul algebra: we take for $A$ the homogeneous
ring of the Veronese curve  $X=\BP^1 \hra \BP^N$. This part contains
no new results; rather it is a review of some classical and modern
theorems related to the Koszul--Chevalley resolution of $A$. This
topic turns out remarkably rich; it revolves around the {\it Gauss
cyclotomic identity}. We see here Euler products, Witt theorem,
Polya theory and a~formula of Polyakov. A particular case of a deep
theorem by Kempf and Bezrukavnikov says that $A$ is a Koszul
algebra. We will see in Theorem~\ref{theorem3.4.4} that the ``numerical'' manifestation of
this fact is precisely the Gauss cyclotomic identity.

\paragraph{1.2.} At the end of this introduction let us say a few words
about the algebras interesting to the physicists; we hope to return
to these questions later on. In his seminal papers N.~Berkovits
considers the algebra of functions on the space of {\it pure
spinors} of dimension $10$, the quotient of the polynomial ring in
$16$ variables $R_0=\BC[\lambda^1,\ldots,\lambda^{16}]$ by the ideal
generated by $10$ quadratic elements
\[
\sum_{\alpha,\beta}  \lambda^\alpha\gamma^m_{\alpha\beta}\lambda^\beta,\qquad
m=0,\ldots, 9.
\]
The corresponding projective variety $i: X = \text{Proj}\, A \hra
\BP^{15}$ is the hermitian symmetric space $X = SO(10,\BC)/U(5)$.
The canonical bundle is $\omega_X = {\cal O} _X(-8)$ where ${\cal
O}_X(1) = i^*{\cal} O_{\BP^{15}}$, cf.~\cite{GKR}.

The Hilbert series of $A$ is equal to
\[
H(A;t) = \frac{1 - 10t^2 + 16t^3 - 16t^5 + 10t^6 - t^8}{(1-t)^{16}} =
\frac{(1+t)(1+4t+t^2)}{(1-t)^{11}}.
\]
The algebra $A$ is Koszul\footnote{The Kempf--Bezrukavnikov theorem says that
the coordinate rings of homogeneous spaces of the form $G/P$ where $G$
is semisimple complex and $P$ is a parabolic,  are Koszul. Another manifestation
of this is the Kapranov's description (\`a la Beilinson and Bernstein--Gelfand--Gelfand) of the derived category of coherent
sheaves on $G/P$.}.

Let  $L = \mathop{\oplus}\limits_{i=1}^\infty  L_i$ be the graded Lie algebra Koszul dual to
$A$ and $C^*(L)$ its Chevalley cochain complex -- the Koszul--Chevalley
resolution of $A$. As a graded algebra it is the tensor product
\[
C^*(L) = S(L[1]^*) = SL^*_1\otimes\Lambda L^*_2\otimes SL^*_3\otimes\cdots
\]
(more precisely the inductive limit of f\/inite products; here $S$
denotes the symmetric algebra and the star denotes the dual space).
Let us consider a graded commutative algebra
\begin{gather}
{\cal C}_1(L) = C^*(L)\otimes \Lambda L^*_1. \label{(1.2.1)}
\end{gather}
One can def\/ine a dif\/ferential on this algebra that will make a dga whose
cohomology will be
\begin{gather}
H^*({\cal C}_1(L)) = {\rm Tor}_*^{R_0}(A,\BC) =: Q .\label{(1.2.2)}
\end{gather}
The dga ${\cal} C_1(L)$ is quasi-isomorphic to the {\it Berkovits
algebra} studied by  Movshev and Schwarz, cf.~\cite{MS}.

The graded commutative algebra $Q$ is called {\it the algebra of
syzygies}. For $i\geq 1$ denote by  $L_{\geq i} \subset L$ the
graded Lie subalgebra $\mathop{\oplus}\limits_{j\geq i} L_j$. It follows from the
above description that $Q$ is isomorphic to the algebra $H^*(L_{\geq
2};\Bbb C)$ -- this is a result contained in \cite{GKR}.

Generalizing this construction, consider for each $i\geq 1$ a graded
commutative algebra
\[
{\cal C}_i(L) = C^*(L)\otimes  \Lambda L^*_1\otimes SL^*_2\otimes
\cdots \otimes F_iL_i^*,
\]
where $F_i = S$ if $i$ is even and  $\Lambda$ otherwise. One expects
that it is possible to introduce a~dif\/ferential on ${\cal C}_i(L)$
which provides a dga whose cohomology is $H^*(L_{\geq i})$. For
example, an algebra quasi-isomorphic to ${\cal C}_3(L)$ has been
studied in a very interesting paper~\cite{AABN}.

The algebra of syzygies $Q$ is a remarkable object. It is f\/inite-dimensional over $\BC$, $Q = \mathop{\oplus}\limits_{i=0}^3 Q_i$ and admits, after~\cite{GKR} a scalar product $Q_i \otimes Q_{3-i}\lra \BC$ compatible
with the multiplication. One can imagine $Q$ as the cohomology of a
smooth compact oriented variety of dimension~$3$.

A few words about the aim of the paper and the novelty of the
results. The physicists are interested in quantum strings
f\/luctuating on the singular space $\text{Spec}(A)$. To pass to
strings one has to study the {\it chiral} analogs of the above
algebras, cf.~\cite{AABN,BD}. The construction of the {\it chiral}
analog of such algebras, the {\it chiralization}, can be
accomplished with the use of free inf\/inite resolutions of the above
algebras. In this paper we put together the properties of Koszul
algebras and their inf\/inite resolutions which contribute to the
properties of the {\it chiral} algebras~\cite{AABN,BD}.  We introduce
and study in Part~I the regularization technique following the
insight from \cite{BN} which allows to def\/ine a part of the above
{\it chiral} algebra structure. The Part~II contains samples of
calculations with the equivariant Hilbert series which one needs to
def\/ine the action of an appropriate Kac--Moody on the {\it
chiralization} of the above algebras. We are planning to return to
{\it chiralization} of these objects in a separate publication.

Otherwise, one can regard this topic from the point of view of the
{\it string topology} of  Chas and Sullivan. The string theory is
the study of spaces of loops. Let $M$ be a closed oriented manifold
of dimension $d$, $\Omega M$ its free loop space. A fundamental
theorem of Chas and Sullivan, cf.~\cite{CS}, says that the homology
of $M$ shifted by $d$
\[
\BH_*(\Omega M) := H_{*+d}(\Omega M)
\]
admits a structure of a {\it Batalin--Vilkovisky} (BV) algebra.
If $M$ is simply connected then $\BH_*(\Omega M)$ is isomorphic
to the Hochschild cohomology
\[
\BH_*(\Omega M) = HH^*(C^*(M),C^*(M)),
\]
where $C^*(M)$ is the complex of singular cochains of $M$, cf.~\cite{CJ,M}.

An algebraic counterpart of the Chas--Sullivan theory is the
following remarkable result (``the cyclic Deligne conjecture''), cf.~\cite{Ka,TZ}: let $A$ be an associative algebra with an
invariant scalar product; then the complex of Hochschild cochains
$CH^*(A)$ admits a structure of a homotopy BV algebra.

Returning to the Berkovits algebras  $\CC_1(L)$, $Q$, cf.~\eqref{(1.2.1)},
\eqref{(1.2.2)}, it would be very interesting to study their Hochschild (as
well as cyclic) cohomology. In view of the previous remarks, it
should be closely related to the state space of the Berkovits
string\footnote{It seems that the Polyakov's gauge f\/ields-strings
correspondence translates in algebra into the assertion: the
Hochschild cohomologies of Koszul dual algebras are isomorphic}.

It is worth mentioning certain analogy between the ``chiral'' and
``topological'' points of view. For example the Deligne conjecture:
the Hochschild cochain complex $CH^*(A)$ of an associative algebra
$A$ is a homotopy Gerstenhaber algebra (or more precisely, an
algebra over the operad $e_2$ of chains of little discs), resembles
the Lian--Zuckerman conjecture~\cite{LZ}: the space of a
topological (i.e.\ $N=2$ supersymmetric) vertex algebra is a homotopy
Gerstenhaber algebra. If an associative algebra is equipped with an
invariant scalar product then $CH^*(A)$ becomes a homotopy BV
algebra. Which complementary structure one needs on a topological
vertex  algebra to become a homotopy BV algebra? In other words,
does there exist a vertex counterpart of the cyclic Deligne
conjecture?

\newpage

\pdfbookmark[1]{First part}{part1}
\section*{First part}

\section{Numerical invariants and regularization}

\begin{flushright}
\begin{minipage}{12cm}
{\it ``Les s\'eries divergentes sont en g\'en\'eral quelque chose de bien fatal,
et c'est une honte qu'on ose y fonder aucune d\'emonstration.
On peut d\'emontrer tout ce qu'on veut en les employant\dots''}
\vspace{2mm}

\qquad \hfill N.-H.~Abel, a letter to Holmboe, January 16, 1826
\end{minipage}
\end{flushright}

\subsection{Semi-free resolutions and Hilbert series}\label{section2.1}

\paragraph{2.1.1.}
Let
\[
A = R_0/(f_1,\ldots,f_M) = \mathop{\oplus}\limits_{i=0}^\infty  A^i,
\]
and $i:  X := \text{Proj}\, A \hra \BP^N := \text{Proj}\, R_0$ be as in the
Introduction; denote by  $\CL = i^*\CO_{\BP^N}(1)$ the corresponding
very ample line bundle.

Set
$h_i = \dim_k  A^i$ and let
\[
H(A;t) = \sum_{i=0}^\infty  h_it^i
\]
be the Hilbert series; we suppose $X$ to be connected, so
$H(A;0) = h_0 = 1$. Let
\begin{gather}
R:\ \cdots \lra R_2 \lra R_1 \lra R_0 \lra A \lra 0 \label{(2.1.1.1)}
\end{gather}
be a semi-free resolution of $A$. For $R_1$ one can take
$R_1 = \mathop{\oplus}\limits_{l=1}^M
R_0\xi_l$, where the variables $\xi_l$ are odd and
$d\xi_l = f_l$, so  $\xi_l \in R_1^{d_i}$.


\begin{Example} Suppose that $f_1, \ldots, f_M$ is a
regular sequence, that is, $X$ is a complete intersection. Then one
can take for $R$ the Koszul complex, $R_i = R_0\otimes \Lambda^i W$,
where $W = \mathop{\oplus}\limits_{l=1}^M k\xi_l$.

The exact sequence \eqref{(2.1.1.1)} immediately gives an expression of the
Hilbert series:
\begin{gather}
H(A;t) = \frac{\prod\limits_{l=1}^M(1 - t^{d_l})}{(1 - t)^{N+1}} =
\frac{\prod\limits_{l=1}^M(1 + t + \cdots + t^{d_l-1})}{(1 - t)^{d+1}} =
\frac{P(t)}{(1 - t)^{d+1}} , \label{(2.1.2.1)}
\end{gather}
where
\begin{gather}
d = N - M  = \dim X. \label{(2.1.2.2)}
\end{gather}
One notes that
\begin{gather}
\deg X = \prod_{l=1}^M\ d_l = P(1). \label{(2.1.2.3)}
\end{gather}
\end{Example}

\paragraph{2.1.2.}
Let us return to the general case. Recall that $R$
as an algebra is a polynomial superalgebra in homogeneous variables.
Let us def\/ine integer numbers
$a_n = $ the number of even generators of polynomial degree $n$
minus the number of odd genera\-tors of polynomial degree $n$.

Then the exact sequence \eqref{(2.1.1.1)} gives a product expression
\begin{gather}
H(t) = H(A;t) = \prod_{n=1}^\infty\left(1 - t^n\right)^{-a_n}
\label{(2.1.3.1)}
\end{gather}
similar to \eqref{(2.1.2.1)}.

This formula shows in particular that the numbers $a_n$ depend only
on $A$ but not on the resolution~$R$.

Similarly to \eqref{(2.1.2.2)} and \eqref{(2.1.2.3)} we want to show that
\begin{gather}
\dim X + 1 = \text{``}a_1 + a_2 + a_3 + \cdots\text{''}, \label{(2.1.3.2)}
\\
\log(\deg X) = - \text{``}(\log 1)\cdot a_1 + (\log 2)\cdot a_2 + (\log
3)\cdot a_3 + \cdots\text{''}. \label{(2.1.3.3)}
\end{gather}
It is natural also to consider the ``higher moments'':
\begin{gather}
\text{``}\sum_{n=1}^\infty\ n^la_n\text{''},\qquad l \geq 1. \label{(2.1.3.4)}
\end{gather}
However, the series at the right hand side are divergent, so one
puts the sums in the quotation marks. Our aim will be to perform the
summation of these series\footnote{The classical sources about
divergent series are the books~ \cite{Ha,R}.}. There are several ways of doing the summation.

For example, one can write a Lambert series
\[
f(t) = \sum_{n=1}^\infty\ a_n\cdot \frac{nt e^{-nt}}{1 - e^{-nt}},
\]
cf.~\cite[App.~IV, (1.1)]{Ha}.

This series diverges if $|t|$ is small; but one can show that if
$|t|$ is suf\/f\/iciently big then the series absolutely converges and
$t^{-1}f(t)$ is a rational function of $y = e^{-t}$. Therefore one
extend~$f(t)$ to the complex plane; this function will be
holomorphic at $t=0$ (i.e.\ and $y=1$). It is natural to def\/ine
\[
\text{``}\sum\ a_n\text{''}:= f(0) = \mathop{\rm res}_{t=0}\frac{f(t)}{t}.
\]
The identity \eqref{(2.1.3.2)} will hold true. Moreover, the higher moments
\eqref{(2.1.3.4)} may be expressed in terms of the coef\/f\/icients of the Taylor
series of $f(t)$ and $0$; this is noted in~\cite{BN}.

Another classical way of regularization is using the Mellin
transform and working with the Dirichlet series. This is what we are
going to do.

\subsection{M\"obius inversion}\label{section2.2}

\paragraph{2.2.1.} Take the logarithm of \eqref{(2.1.3.1)}:
\begin{gather}
\log\ H(t) = - \sum_{n=1}^\infty  a_n\log\left(1 - t^n\right) =
\sum_{n=1}^\infty \sum_{l=1}^\infty a_n\frac{t^{ln}}{k} =
 \sum_{m=1}^\infty\left(\sum_{l|m} \frac{a_l}{m/l}\right)\cdot
t^m \label{(2.2.1.1)}
\end{gather}
and then the derivative:
\begin{gather}
\frac{tH'(t)}{H(t)} = \sum_{n=1}^\infty \frac{na_n t^n}{1 - t^n} =
 \sum_{m=1}^\infty\left(\sum_{l|m} la_l\right)\cdot t^m.
\label{(2.2.1.2)}
\end{gather}
In other words, if one denotes
\begin{gather}
\frac{tH'(t)}{H(t)} = \sum_{m=1}^\infty b_m t^m, \label{(2.2.1.3)}
\end{gather}
then
\begin{gather}
b_m = \sum_{l|m} la_l. \label{(2.2.1.4)}
\end{gather}

\paragraph{2.2.2.} Recall that the M\"obius function $\mu:  \BN_+ =\{1,
2, \ldots \} \lra \{-1, 0, 1\}$ is def\/ined by: $\mu(1) = 1$, $\mu(n) =
(-1)^l$ if $n = p_1p_2\cdots p_l$ is a product of $l$ distinct prime
numbers, and $\mu(n) = 0$ if $n$ contains squares.

Another def\/inition is by a generating series: if one def\/ines
the Riemann  $\zeta$ function by the Euler product
\[
\zeta(s) = \sum_{n=1}^\infty n^{-s} = \prod_{p\ \text{prime}} \left(1 -
p^{-s}\right)^{-1},
\]
then
\begin{gather}
\zeta(s)^{-1} = \prod_{p\ \text{prime}} \left(1 - p^{-s}\right) =
\sum_{n=1}^\infty \mu(n) n^{-s}. \label{(2.2.2.1)}
\end{gather}

{\it The M\"obius inversion formula} says that if
$f:  \BN_+ \lra \BC$ is a function and if a function
$g:  \BN_+ \lra \BC$ is def\/ined by
\[
g(m) = \sum_{l|m} f(l),
\]
then
\begin{gather}
f(m) = \sum_{l|m} \mu(m/l)g(l). \label{(2.2.2.2)}
\end{gather}
Applying this formula to the function $f(m) = ma_m$, we get
\begin{gather}
ma_m = \sum_{l|m} \mu(m/l)b_l. \label{(2.2.2.3)}
\end{gather}

\subsection{Dirichlet series}\label{section2.3}

\paragraph{2.3.1.} Let $P(t)$ be a polynomial with complex coef\/f\/icients such
that $P(0) = 1$, so one can write
\begin{gather}
P(t) = \prod_{l=1}^p (1 - \alpha_lt). \label{(2.3.1.1)}
\end{gather}
Consider the product \eqref{(2.1.3.1)} for $P(t)$:
\[
P(t) = \prod_{n=1}^\infty \left(1-t^n\right)^{-a_n}.
\]
So if
\[
\frac{tP'(t)}{P(t)} = \sum_{m=1}^\infty b_mt^m,
\]
then
\[
m a_m = \sum_{l|m} \mu(m/l)b_l.
\]
On the other hand
\begin{gather}
\frac{tP'(t)}{P(t)} = - \sum_{r=1}^p \frac{\alpha_rt}{1 - \alpha_rt} =
 - \sum_{r=1}^p \sum_{m=1}^\infty \alpha_r^mt^m,
\end{gather}
i.e.
\[
b_m = - \sum_{r=1}^p \alpha_r^m,
\]
wherefrom
\[
a_m = - m^{-1}\sum_{r=1}^p\sum_{l|m} \mu(m/l)\alpha_r^l.
\]
One puts
\[
\rho(P) =  \max_r |\alpha_r|.
\]
It follows that if $\rho(P) > 1$, then $|a_m|$ grows as fast as
$m^{-1}\rho(P)^m$.

\paragraph{2.3.2.} Let us consider the Dirichlet series
\[
\sum_{n=1}^\infty a_n n^{-s}.
\]
We see that if $\rho(P) \leq 1$ then this series absolutely
converges for ${\rm{\fR}e}(s) > 1$, and if $\rho(P) > 1$, it diverges for
all $s$; in this case (which in fact is interesting to us) one has
to do something else.

Let us write formally with \cite{BN}:
\begin{gather*}
- \sum_{n=1}^\infty a_n n^{-s} = \sum_{n=1}^\infty \sum_{r=1}^p\sum_{l|n} \mu(n/l)\alpha_r^l n^{-1-s} =
 \sum_{r=1}^p \sum_{l, m = 1}^\infty \mu(m)\alpha_r^l (lm)^{-1-s} \\
 \phantom{- \sum_{n=1}^\infty a_n n^{-s}}{} =
\sum_{m=1}^\infty \mu(m) m^{-1-s}\cdot
\sum_{r=1}^p \sum_{l = 1}^\infty \alpha_r^l l^{-1-s} =
 \zeta(s+1)^{-1}\sum_{r=1}^p\sum_{l = 1}^\infty \alpha_r^l l^{-1-s}.
\end{gather*}
Now if we rewrite, after Riemann,
\[
l^{-1-s} = \frac{1}{\Gamma(s+1)}\int_0^\infty e^{-lt}t^{s}dt,
\]
wherefrom
\[
\sum_{l = 1}^\infty \alpha_r^l l^{-1-s} =
\frac{1}{\Gamma(s+1)}\int_0^\infty
\sum_{l = 1}^\infty \alpha_r^le^{-lt}t^{s}dt =
 \frac{1}{\Gamma(s+1)}\int_0^\infty  \frac{\alpha_re^{-t}t^{s}}
{1 - \alpha_re^{-t}}dt
\]
(justif\/ied if $|\alpha_r| < 1$). Because
\[
\frac{e^{-t}}{1 - \alpha_re^{-t}} =
\frac{1}{e^t - \alpha_r},
\]
the integral
\[
\int_0^\infty \frac{\alpha_re^{-t}t^{s}}
{1 - \alpha_re^{-t}}dt
\]
absolutely converges for ${\rm{\fR}e}(s) > - 1$ if $\alpha_r \not\in
[0,1]$ (since $P(1) = 1$, one has $\alpha_r \neq 0$ automatically).
One notes that
\begin{gather}
\sum_{r=1}^p  \frac{\alpha_re^{-t}} {1 - \alpha_re^{-t}} = -
\frac{e^{-t}P'(e^{-t})}{P(e^{-t})}. \label{(2.3.2.1)}
\end{gather}
On the other hand, let us write the functional equation for
$\zeta(s)$ (cf.~\cite[13.151]{WW}):
\[
2^{1-s}\Gamma(s)\zeta(s)\cos(\pi s/2) =
\pi^s\zeta(1-s)
\]
and replace $s$ by $s+1$:
\[
2^{-s}\Gamma(s+1)\zeta(s+1)\cos(\pi (s+1)/2) =
\pi^{s+1}\zeta(-s),
\]
i.e.
\[
\frac{1}{\Gamma(s+1)\zeta(s+1)} = - 2^{-s}\pi^{-s-1}\sin(\pi s/2)\zeta(-s)^{-1}.
\]
It follows:
\begin{gather*}
\sum_{n=1}^\infty a_n n^{-s} =  \frac{1}{\Gamma(s+1)\zeta(s+1)}
\int_0^\infty \frac{e^{-t}P'(e^{-t})}{P(e^{-t})} t^s dt\\
\phantom{\sum_{n=1}^\infty a_n n^{-s}}{}
= - \frac{2^{-s}\pi^{-s-1}\sin(\pi s/2)}{\zeta(-s)}
\int_0^\infty \frac{e^{-t}P'(e^{-t})}{P(e^{-t})} t^s dt.
\end{gather*}

\paragraph{2.3.3.} Consider the last integral
\[
I_P(s) := \int_0^\infty \frac{e^{-t}P'(e^{-t})}{P(e^{-t})} t^s dt.
\]
It is well def\/ined and represents a holomorphic function of $s$ on
the half-plane ${\rm{\fR}e}(s) > - 1$ as soon as the condition \eqref{$*$} below is
verif\/ied:
\begin{gather}
\mbox{The roots of $P(t)$ do not lie in the segment $[0,1]$.} \label{$*$}
\end{gather}
As it was kindly pointed out by one of the referees the condition
\eqref{$*$} holds. Indeed, since the algebra $A$ is af\/f\/ine, it has
polynomial growth, so that the radius of convergence of the series
$H(A,t)$ is  $0\leq t< 1$. Because the coef\/f\/icients of the formal
power series $H(A,t)$ are nonnegative, we  have $H(A, t )
>0$ for all $0 \leq t < 1$. It follows that $P(A,t) > 0$ with $H(A,
t )$. Also, the is function $H(A,t)$ is rational (Hilbert--Serre), so
that it has a pole at $t=1$ of the order $d+1$. This means that the
value $P(A,1)$ is nonzero as well.

\paragraph{2.3.4.} Following Riemann and
using the notations of \cite{WW}, consider the integral
\begin{gather}
J_P(s) = \int_\infty^{(0+)}  \frac{e^{-t}P'(e^{-t})}{P(e^{-t})}
(-t)^s dt, \label{(2.3.4.1)}
\end{gather}
where $\int_\infty^{(0+)}$ denotes the integral along the contour
\begin{gather}
C = \{ t = a + \epsilon i,\, \infty > a \geq \epsilon\} \cup
\{ t = \sqrt 2\epsilon e^{i\theta},\, \pi/4 \leq \theta \leq 7\pi/4\}\nonumber\\
\phantom{C =}{}
\cup \{ t = a - \epsilon i,\, \epsilon \leq a < \infty\},
\label{(2.3.4.2)}
\end{gather}
where $\epsilon > 0$ is suf\/f\/iciently small. Here $(-t)^s =
e^{s\log(-t)}$, and we use the branch of the logarithm $\pi \geq
\arg t \geq -\pi$ on the small circle. Then
\[
J_P(s) = 2i\sin (\pi s) I_P(s),
\]
cf.~\cite[12.22]{WW}. But $J_P(s)$ is an entire function.
Therefore, if one sets
\begin{gather}
z_P(s) = - \frac{2^{-s}\pi^{-s-1}\sin(\pi s/2)}{\zeta(-s)}\cdot
\frac{1}{2i\sin(\pi s)} J_P(s) =
 \frac{i\cdot 2^{-s-2}\pi^{-s-1}}{\cos(\pi s/2)\zeta(-s)} J_P(s),
\label{(2.3.4.3)}
\end{gather}
this def\/ines the analytic continuation of \eqref{(2.3.4.1)} to a meromorphic
function on the complex plane.

\paragraph{2.3.5.} Let us  return to the situation of Section~\ref{section2.1}; it is known that
the Hilbert series is a~rational function of the form
\[
H(A;t) = \frac{P(A;t)}{(1 - t)^{d+1}},
\]
where $d = \dim X = \dim A - 1$ and $P(A;t)$ is a polynomial with
integer coef\/f\/icients (cf.\ a simple Example~\ref{example2.3.7}). One has
$P(A;0) = H(A;0) = 1$ since we  suppose that $X$ connected.

In view of the preceding discussion we def\/ine
\begin{gather}
z(X,\CL;s) = \frac{1}{\Gamma(s+1)\zeta(s+1)}\int_0^\infty
\frac{e^{-t}H'(A;e^{-t})}{H(A;e^{-t})}t^s dt  =
 d + 1 + z_P(s) \label{(2.3.5.1)}
\end{gather}
(because $\Gamma(s+1)\zeta(s+1) = \int_0^\infty  (t^s/(e^t - 1)) dt$).
So if
\begin{gather}
H(A;t) = \prod_{n=1}^\infty  (1 - t^n)^{-A_n}, \label{(2.3.5.2)}
\end{gather}
then
\[
z(X,\CL;s) = \text{``}\sum_{n=1}^\infty A_n n^{-s}\text{''},
\]
the quotation marks mean that the sum is regularized. We def\/ine ``the
Dirichlet summation'':
\[
\text{``}\sum_{n=1}^\infty  A_n\text{''}:= z(X,\CL;0),
\]
and
\[
\text{``}\sum_{n=1}^\infty  \log n\cdot A_n\text{''}:= - z'(X,\CL;0),
\]

\paragraph{2.3.6.} We have $\zeta(0) = - 1/2$ and $\zeta'(0) = -
\log\sqrt{2\pi} \neq 0$ (cf.~\cite[Ch.~VII, \S~9]{Weil}). The integral
$I_P(0)$ is well def\/ined (below we shall compute its value), so the
formula \eqref{(2.3.5.1)} gives $z_P(0) = 0$, wherefrom
\[
z(X,\CL;0) = d + 1,
\]
cf.~\eqref{(2.1.3.2)}. On the other hand, $z'(X,\CL;s) = z'_P(s)$, and a
unique term giving a nontrivial contribution into $z'(X,\CL;0)$, is
\[
- \left(\frac{2^{-s}\pi^{-s-1}\cdot (\pi/2)\cos(\pi s/2)}{\zeta(-s)}
\right)\biggr|_{s=0}\cdot I_P(0).
\]
The f\/irst factor gives $1$, whereas
\begin{gather*}
I_P(0) = \int_0^\infty \frac{e^{-t}P'(e^{-t})}{P(e^{-t})}  dt =
- \int_0^\infty \frac{d\log P(e^{-t})}{dt} dt
\\
\phantom{I_P(0)}{} = - \log P(e^{-t}) \bigl|_0^\infty = - \log P(0) + \log P(1) = \log P(1),
\end{gather*}
wherefrom
\[
z'(X,\CL;0) = \log P(1).
\]
But it is known that
$P(1) = \deg X$, which implies
\[
\deg X = e^{z'(X,\CL;0)},
\]
cf.\ \eqref{(2.1.3.3)}.


\begin{Example}[Veronese curve and the Weil zeta
function]\label{example2.3.7}
Let $X = \BP^1$, $\CL = \CO(q+1)$, $q \geq 1$. Then $A^n =
\Gamma(X,\CO((q+1)n)$, which implies
\[
H(A;t) = \sum_{n=0}^\infty ((q+1)n + 1)t^n = \frac{1 + qt}{(1 - t)^2},
\]
and
\[
P(t) = 1 + qt = \prod_{m=1}^\infty (1 - t^m)^{p_m(-q)},
\]
where
\[
p_m(x) = \frac{1}{m}\sum_{l|m} \mu(m/l)x^l.
\]
It is well known that if $q$ is a prime power then
$p_m(q)$ is equal to the number of monic irreducible polynomials of degree $m$
in $\BF_q[T]$.
In fact, we have
\[
Z_{\BF_q}(\BA^1,t) = \frac{1}{P(-t)} = \prod_{m=1}^\infty (1 - t^m)^{-p_m(q)}
\]
-- this is the Euler product of the zeta function of the af\/f\/ine line
over the f\/inite f\/ield $\BF_q$.
\end{Example}

(The idea that the Hilbert series of a variety is equal to the Weil
zeta function of another one seems very strange. Cf.~\cite{Gol}
however\footnote{We are grateful to Yu.I.~Manin who showed us this
very interesting paper.}.)

\subsection{Values at negative points}\label{section2.4}

\paragraph{2.4.1.} Let $s = -m$ be a negative integer. In this case one can
close the contour in the integral~\eqref{(2.3.4.1)}:
\[
J_P(-m) = (-1)^m\int_{|t| = \epsilon}  Q(t)t^{-m}dt =
2\pi i\cdot (-1)^m\mathop{\rm res}\limits_{t=0} Q(t)t^{-m},
\]
where
\[
Q(t) = \frac{e^{-t}P'(e^{-t})}{P(e^{-t})} = - \frac {d\log P(e^{-t})}{dt}.
\]
So if
$Q(t) = q_0 + q_1t/1! + q_2t^2/2! + \cdots$
is the Taylor series at $0$ then
\[
J_P(-m-1) = 2\pi i\cdot (-1)^{m+1}\frac{q_m}{m!}.
\]
Here are some f\/irst values:
\begin{gather*}
q_0 = \frac{P'(1)}{P(1)},\qquad q_1 = - q_0 + q_0^2 - \frac{P''(1)}{P(1)},
\\
q_2 = q_0 - 3q_0^2 + 2q_0^3
+ 3(1 - q_0)\frac{P''(1)}{P(1)} + \frac{P'''(1)}{P(1)}.
\end{gather*}

\paragraph{2.4.2.} One can say this in a dif\/ferent way: if
\[
P(t) = 1 + c_1t + \cdots + c_Dt^D,
\]
and we imagine the numbers $c_i$ as ``Chern classes''
then $q_m$ will be the coef\/f\/icients of the logarithm of the Todd genus \dots

\paragraph{2.4.3.} Now applying \eqref{(2.3.4.3)}: for $m = 1$, the function $\zeta(s)$
has a simple pole at $s = 1$ with the residue $1$, whence
\begin{gather}
z_P(-1) = \frac{2P'(1)}{P(1)}. \label{(2.4.3.1)}
\end{gather}


\begin{Remark} Suppose that $P(t)$ is a ``reciprocal''
polynomial, i.e.
\[
t^{\deg P}P(1/t) = P(t),
\]
which happens rather often. In this case it is easy to see that
\[
\frac{2P'(1)}{P(1)} = \deg P,
\]
so $z_P(-1) = \deg P$.
\end{Remark}

\paragraph{2.4.4.}
 For $m > 1$ one has to distinguish two cases.

(a) If $m = 2l$ is even then
$\cos(\pi l) = (-1)^l$, and after Euler,
\[
\zeta(2l) = (-1)^{l-1}\frac{(2\pi)^{2l}}{2(2l)!} b_{2l},
\]
where the Bernoulli numbers are def\/ined by the generating series
\[
\frac{S}{e^S - 1} =  1 - \frac{S}{2} + \sum_{l=1}^\infty
\frac{b_{2l}}{(2l)!}S^{2l}.
\]
Here are some f\/irst values:
\[
b_2 = \frac{1}{6},\qquad  b_4 = - \frac{1}{30},\qquad
b_6 = \frac{1}{42},\qquad  b_8 = - \frac{1}{30}.
\]
It follows:
\[
z_P(-2l) = \frac{2l q_{2l-1}}{b_{2l}}.
\]
So if the numbers  $a_n$ are def\/ined by \eqref{(2.3.5.2)} then
\begin{gather}
\text{``}\sum_{n=1}^\infty n^{2l}a_n\text{''} = z(X,\CL;-2l) = d + 1 +
\frac{2l q_{2l-1}}{b_{2l}}
\label{(4.5.1)}
\end{gather}
-- a formula found empirically in \cite[(4.21)]{BN}.

(b) If $m = 2l + 1$, $l > 0$ then $z_P(s)$ has a simple pole at $s = - m$, with the
residue
\[
\mathop{\rm res}\limits_{s = -2l-1} z_P(s) = (-1)^{l}\frac{(2\pi)^{2l}q_{2l}}
{(2l+1)!\zeta(2l+1)}.
\]

\paragraph{2.4.5.}
 {\it Hurwitz $zeta$ function.}

(a) We have $J_P(s) = \sum\limits_{r=1}^p J_{\alpha_r}(s)$ where
\[
J_\alpha(s) = - \int_\infty^{(0+)} \frac{\alpha}{e^t - \alpha} (-t)^s dt.
\]
Let us choose a value $\beta = \log  \alpha$. The expression under the integral
$f_\alpha(t,s) = (\alpha/e^t - \alpha) \cdot (-t)^s$
has the poles at points $t_n = \beta + 2\pi i n$, $n\in \BZ$. By the
Cauchy formula,
\[
J_\alpha(s) + 2\pi i\sum_{n= -\infty}^\infty  \mathop{\rm res}\limits_{t=t_n} f_\alpha(t,s)
= 0
\]
if ${\rm{\fR}e}(s) < -1$, then
\begin{gather}
J_\alpha(s) = 2\pi i\sum_{n= -\infty}^\infty  (-\beta - 2\pi i n)^s,
\label{(2.4.6.1)}
\end{gather}
which gives $J_P(s)$ in the form of a series absolutely convergent
for ${\rm{\fR}e}(s) < -1$, expressible in terms of the Hurwitz $zeta$
functions.

{\samepage (b) \begin{Example} Take $P(t) = 1 + pt$, where $p\in \BZ_{>0}$, cf.\ Example~\ref{example2.3.7}. Then $q_0 = p/(p+1)$, $q_1 = - p/(p+1)^2$ (cf.~Section~2.4.1), so
\[
\frac{J_P(-2)}{2\pi i} = - \frac{p}{(p+1)^2}.
\]
On the other hand, take
$\beta = \log (-p) = \log p + \pi i$, and the above series will be written as
\[
\sum_{n=-\infty}^\infty (-\log p - \pi i - 2\pi i n)^{-2} =
\sum_{n=0}^\infty \left\{(-\log p + (2n+1)\pi i)^{-2} +
(-\log p - (2n+1)\pi i)^{-2}\right\},
\]
so one came to an identity
\begin{gather}
\sum_{n=0}^\infty \left\{(-\log p + (2n+1)\pi i)^{-2} + (-\log p -
(2n+1)\pi i)^{-2}\right\} = - \frac{p}{(p+1)^2}. \label{(2.4.6.2)}
\end{gather}
For example, for
$p=3$, $-p/(p+1)^2 = - 0,1875$, whereas an approximate value
given by MAPLE is
\[
\sum_{n=0}^{1000} \left\{(-\log p + (2n+1)\pi i)^{-2} +
(-\log p - (2n+1)\pi i)^{-2}\right\} = - 0,187449\dots;
\]
we see that the convergence is very slow.
\end{Example}}

(c) If $s = -1$, the series \eqref{(2.4.6.1)} in our example is still ``Eisenstein summable'',
and on gets a~rational value:
\begin{gather}
\sum_{n=0}^\infty \left\{(-\log p + (2n+1)\pi i)^{-1} + (-\log p -
(2n+1)\pi i)^{-1}\right\} = - \frac{2(p-1)}{p+1}. \label{(2.4.6.3)}
\end{gather}
This is an easy corollary of the decomposition of $\cot z$ into
simple fractions\footnote{We thank Oleg Ogievietsky who has noted
this.}. But this value is dif\/ferent from  $- q_0 = - p/(p+1)$: we
have ``an anomaly''.

\paragraph{2.4.6.}
Suppose that our variety $X$ is smooth, $H^i(X,\CL^{\otimes
n}) = 0$ for all $i > 0$ and $n \geq 0$. Then $A^n =
\Gamma(X,\CL^{\otimes n})$, and we can switch on the
Riemann--Roch--Hirzebruch~\cite{H1}:
\[
h_n = \dim H^0(X,\CL^{\otimes n}) = \int_X e^{nc_1(\CL)}{\rm Td}(\CT_X),
\]
where $\CT_X$ denotes the tangent bundle,
${\rm Td}$ the Todd genus, given for line bundle
$E$ by the formula,
\[
{\rm Td}\,(E) = \frac{c_1(E)}{1 - e^{-c_1(E)}},
\]
so  ${\rm Td}\,(\CT_X) = 1 + {\rm Td}_1(\CT_X) + {\rm Td}_2(\CT_X) + \cdots$, where
\[
{\rm Td}_1(\CT_X) = \frac{c_1(\CT_X)}{2},\qquad
{\rm Td}_2(\CT_X) = \frac{c_1^2(\CT_X) + c_2(\CT_X)}{12}.
\]
One writes:
$e^{nc_1(\CL)} = \sum\limits_{i=0}^\infty c_1(\CL)^i n^i/i!$,
so $h_n = R(n)$, where
\begin{gather*}
R(t) = \biggl(\int_X c_1(\CL)^d\biggr)\cdot \frac{t^d}{d!}  +
\biggl(\int_X c_1(\CL)^{d-1}c_1(\CT_X)\biggr)\cdot \frac{t^{d-1}}{2(d-1)!} +
\\
\phantom{R(t) =}{} + \biggl(\int_X c_1(\CL)^{d-2}(c_1(\CT_X)^2 + c_2(\CT_X))\biggr)\!
\cdot \frac{t^{d-2}}{12(d-2)!} + \cdots =
r_d\frac{t^d}{d!} + r_{d-1}\frac{t^{d-1}}{(d-1)!} + \cdots
\end{gather*}
is a polynomial of degree
$d = \dim X$, the Hilbert polynomial of the ring $A$; it is a polynomial with
rational coef\/f\/icients which takes integer values at integer argument, whence
\[
R(t) = \sum_{i=0}^d  (-1)^{d-i}e_{d-i}\binom{t+i}{i},
\]
where $e_0, e_1, \ldots, e_d \in \BZ$. For example,
\begin{gather*}
e_0 = r_d = \int_X c_1(\CL)^d,
\\
e_1 = \frac{d+1}{2} r_d - r_{d-1} = \frac{1}{2}\int_X \bigl((d+1)c_1(\CL)^d - c_1(\CL)^{d-1}c_1(\CT_X)\bigr),
\\
e_2 = r_{d-2} - \frac{d}{2}r_{d-1} + \frac{(d+1)(3d-2)}{24}r_d
\end{gather*}
(we use $\sum\limits_{1\leq i <j \leq d} ij = d(d^2-1)(3d+2)/24$).

On the other hand it is known that
\[
e_i = \frac{P^{(i)}(1)}{i!}.
\]
It follows:
\[
e^{z'(X,\CL;0)} = P(1) = e_0 = \int_X c_1(\CL)^d,
\]
so one f\/inds the degree of $X$. Next
\begin{gather*}
z_P(-1) = - \frac{2P'(1)}{P(1)} = - \frac{2 e_1}{e_0} =
- \frac{1}{e_0}\int_X
\bigl((d+1)c_1(\CL)^d - c_1(\CL)^{d-1}c_1(\CT_X)\bigr)
\\
\phantom{z_P(-1)}{}
= - d - 1 + \frac{1}{e_0}\int_X  c_1(\CL)^{d-1}c_1(\CT_X),
\end{gather*}
i.e.
\[
z(X,\CL;-1) = \frac{1}{e_0}\int_X  c_1(\CL)^{d-1}c_1(\CT_X),
\]
etc.

One can put these calculations into a formal ``generating function''
for the moments
\[
 M_l=\text{``}\sum_{n=1}^\infty n^{l}a_n\text{''}.
 \]
Note that from the formulas above and using the expansion
\[
\log\left(1-e^t\right)=\log(-t)+\frac{t}{2}+\sum_{l=1}^{\infty} \frac{b_{2l}}{2l(2l)!} t^{2l}
\]
one concludes the formal identity:
\begin{gather}
-\log P(e^t)=\text{``} \sum_{l=0}^{\infty} \log(l)a_l+M_1\frac{t}{2}
+\sum_{l=1}^{\infty} (M_{2l})\frac{b_{2l}}{2l(2l)!} t^{2l}\text{''}.\label{(2.4.7.1)}
\end{gather}

The following elegant formula is due to F.~Hirzebruch \cite{H2}: the
polynomial $P(t)$ is a characteristic number def\/ined as
\[
P(t)=\int_X \frac{(1-t)e^{(1-t)c_1}\widehat {A}((1-t)^2p_1,\dots,(1-t)^{2i}p_{2i},\dots)}{1-te^{(1-t)g}},
\]
where $c_1$, $p_1$, $p_2$, $\dots$ are the f\/irst Chern class and the
Pontryagin classes of the manifold~$X$.

From this it is immediate to express the left hand side of \eqref{(2.4.7.1)}
as a characteristic number. Indeed
\[
 \log P(e^t)=\log\int_X \frac{(1-e^t)e^{(1-e^t)c_1}\widehat {A}((1-e^t)^2p_1,\dots,(1-e^t)^{2i}p_{2i},\dots)}{1-e^{t+(1-e^t)}}.
\]
In order to expand the right hand side into a series of $t$ observe
that the constant term of the series under the logarithm is
\[
\int_X \frac{1}{1-g}=P(1).
\]
Therefore
\[
\log\frac{1}{\int_X \frac{1}{1-g}}\int_X \frac{(1-e^t)e^{(1-e^t)c_1}\widehat {A}((1-e^t)^2p_1,\dots,(1-e^t)^{2i}p_{2i},\dots)}{1-e^{t+(1-e^t)}}
\]
is a series in $t$ equal to $\log P(e^t)-\log P(1)$. This gives the
desired generating function.


 \begin{Example} Let $X_g$ be the moduli space of
semi-stable vector bundles of rank $2$ with tri\-vial determinant over
a Riemann surface of genus $g$; it carries the canonical {\it
determinant} line bundle~$\CL_g$, cf.~\cite{BL}. Consider a graded
algebra
\[
A_g = \mathop{\oplus}\limits_{i=0}^\infty H^0\big(X_g,\CL_g^{\otimes i}\big)
\]
(we thank Peter Zograf who proposed to consider this example).

The coef\/f\/icients of the Hilbert series $H_g(t) = H(A_g;t)$ can be
calculated using the Verlinde formula~\cite{BL}. Here are some
f\/irst examples, cf.~\cite{Z}:
\begin{gather}
H_2(t) = \frac{1}{(1-t)^4},\qquad H_3(t) = \frac{1+t+t^2+t^3}{(1-t)^7} =
\frac{1-t^4}{(1-t)^8}, \label{(2.4.8.1)}
\\
H_4(t) = \frac{1 + 6t + 21t^2 + 40t^3 + 21t^4 + 6t^5 +
t^6}{(1-t)^{10}}. \label{(2.4.8.2)}
\end{gather}
The highest coef\/f\/icient
$r_d/d!$ of the Hilbert polynomial has been calculated in
[Wi], who found the value $2\zeta(2g-2)/(2\pi^2)^{g-1}$; here
$d = \dim X_g = 3g-3$, whence
\begin{gather}
P_g(1) = (-1)^g \frac{(3g-3)!}{(2g - 2)!}b_{2g-2} = (-1)^{g+1}
\frac{(3g-3)!}{(2g - 1)!}\zeta(-2g+3) \label{(2.4.8.3)}
\end{gather}
(one denotes $P_g(t) = H_g(t)(1-t)^{d+1})$. For example, for $g=4$
one f\/inds $P_4(1) = 96$, which is compatible with~\eqref{(2.4.8.2)}.

The above discussion gives a strange expression of this number as a
regularized inf\/inite product (starting from $g=4$), of the form
$\prod\limits_{n=1}^\infty n^{a_n}$. For example, the beginning of the
product for $P_4(t)$ will be:
\begin{gather}
P_4(t) = (1 -
t)^{-6}\big(1-t^3\big)^{16}\big(1-t^4\big)^{9}\big(1-t^5\big)^{-144}\big(1-t^6\big)^{360}\nonumber\\
\phantom{P_4(t) =}{}\times
\big(1-t^8\big)^{-2259}\big(1-t^9\big)^{3920}\cdots, \label{(2.4.8.4)}
\end{gather}
whence
\[
\zeta(-5) = - \frac{9!}{7!}\cdot \text{``}1^{6}3^{-16}4^{-9}5^{144}6^{-360}
8^{2259}9^{-3920}\cdots\text{''}.
\]
\end{Example}

\pdfbookmark[1]{Second part}{part2}
\section*{Second part}

\section[Koszulness and infinite resolutions]{Koszulness and inf\/inite resolutions}\label{section3}

\subsection{The Veronese ring}\label{section3.1}

\paragraph{3.1.1.} Let us f\/ix an integer $b \geq 0$, and consider the Veronese
embedding
\[
i_b:\ \BP^1 \lra \BP^{b+1}
\]
def\/ined in coordinates by the formula
\[
i_b(u_0:u_1) = (x_0:\cdots :x_{b+1}) := \big(u_0^{b+1}:u_0^bu_1:\cdots
:u_1^{b+1}\big).
\]
If $t = u_1/u_0 \in \BA^1 \subset \BP^1$, then
\[
i_b(t) = \big(t,t^2,\ldots,t^{b+1}\big) \in \BA^{b+1} \subset \BP^{b+1}.
\]
The image $X_b:=i_b(\BP^1) \subset \BP^{b+1}$ of this embedding is
called the Veronese curve (or the moment curve, in view of the last formula).

This curve can be def\/ined in  $\BP^{b+1}$ by the equations $x_ix_j -
x_kx_l = 0$ if $i+j = k+l$. A minimal system of equations consists
of $b(b+1)/2$ quadratic equations:
\[
f_{ij}:= x_ix_j - x_{i+1}x_{j-1} = 0,\qquad 0\leq i \leq b-1,\qquad i+2 \leq j
\leq b+1.
\]

 \begin{Example} $i_0 = \text{Id}_{\BP^1}$. If $b=1$, $X_1
\subset \BP^2= \{(x_0:x_1:x_2)\}$ is def\/ined by one equation:
\[
x_0x_2 - x_1^2 = 0.
\]
The curve $X_3 \subset \BP^4$ is def\/ined by $3$ equations
\[
x_0x_2 - x_1^2 = 0;\qquad x_0x_3 - x_1x_2 = 0;\qquad x_1x_3 - x_2^2 = 0.
\]
\end{Example}


\paragraph{3.1.2.}
We have $i_b^*\CO_{\BP^{b+1}}(1) = \CO_{\BP^1}(b+1)$, all
the higher cohomology  of this sheaf vanish, and $X_b = \text{Spec\
Proj}\, A_b$, where
\[
A_b = \mathop{\oplus}\limits_{n=0}^\infty \Gamma\big(\BP^1,\CO((b+1)n\big) =
\BC[x_0,\ldots,x_{b+1}]/(f_{ij}).
\]
The Hilbert series of the above algebra is as follows:
\begin{gather}
H(A_b;t) = \sum_{n=0}^\infty ((b+1)n + 1)t^n = \frac{1 +
bt}{(1-t)^2}. \label{(3.1.3.1)}
\end{gather}
It is not dif\/f\/icult to f\/igure out the  $q$-analogue of \eqref{(3.1.3.1)}: if
we set $[b]_q := (q^b - q^{-b})/(q - q^{-1})$, then:
\begin{gather}
H_q(A_b;t):= \sum_{n=0}^\infty [(b+1)n + 1]_qt^n = \big(q -
q^{-1}\big)^{-1}\cdot \sum_n \big\{q q^{(b+1)n}t^n - q^{-1}q^{-(b+1)n}t^n\big\}
\nonumber\\
\phantom{H_q(A_b;t)}{}
= \big(q - q^{-1}\big)^{-1}\cdot\left(\frac{q}{1 - q^{b+1}t} -
\frac{q^{-1}}{1 - q^{-b-1}t}\right) = \frac{1 + [b]_qt}{(1 -
q^{b+1}t)(1 - q^{-b-1}t)}. \label{(3.1.3.1)_q}
\end{gather}
We will use these formulas later in Section~3.4.4.


\begin{Example} \label{example3.1.4} For $b=1$ we have,
\begin{gather}
H(A_1;t) = \frac{1+t}{(1-t)^2} = \frac{1-t^2}{(1-t)^3},
\label{(3.1.4.1)}
\end{gather}
the ring $A_1$ is a cone:
\[
A_1 = \BC[x_0,x_1,x_2]/\big(x_0x_2 - x_1^2\big) = B/(f).
\]
It admits a dga resolution of length $1$: the Koszul complex:
\[
K_\cdot (B;f):\ 0 \lra B\cdot e \lra B \lra A \lra 0,\qquad d(e) = f.
\]
Therefore
\[
K_\cdot (B;f) = B\otimes \Lambda^\cdot \langle e\rangle = \BC[x_0,x_1,x_2;e]
\]
as a graded algebra; the {\it homological} degree of $x_i$ is $0$
and of $e$ is $1$. These correspond to the exponents $-3$ and $1$ in
\eqref{(3.1.4.1)}.
\end{Example}

There is another interpretation of the Koszul complex. Let $V$ be a
vector space $\mathop{\oplus}\limits_{i=0}^2  \BC\cdot x_i$, then $B = S^\cdot V$.
Consider the dual space  $V^* = \mathop{\oplus}\limits_{i=0}^2 \BC\cdot y_i$. Def\/ine
a graded Lie algebra $L$ on $3$ generators $y_i$, $i = 0, 1, 2$, of
degrees $1$, obeying $5$ relations:
\begin{gather*}
[y_0,y_0] = 0,\qquad [y_0,y_1] + [y_1,y_0] = 0,\qquad [y_0,y_2] + [y_1,y_1] +
[y_2,y_0] = 0,\\
[y_1,y_2] + [y_2,y_1] = 0,\qquad [y_2,y_2] = 0.
\end{gather*}
It is concentrated in degrees $1$ and $2$: $L = L_1 \oplus L_2$,
where $L_1 = V^*$ and $L_2 = \BC\cdot [y_1,y_1]$.

If $L' = \BC\cdot y_1 \oplus \BC\cdot [y_1,y_1] \subset L$ is a Lie
subalgebra generated by $y_1$, then $L'$ is free (sic!); it is an
ideal, and the quotient algebra  $\bar L = L/L'$ is Abelian, on $2$
generators $\bar y_0$, $\bar y_2$.

Consider the complex of the Chevalley cochains:
\[
C^\cdot (L) = \Lambda^\cdot (L^*) = S^\cdot L^*_1 \otimes \Lambda^\cdot L^*_2 = S^\cdot V
\otimes \Lambda^\cdot \langle e\rangle,\qquad  e = [y_1,y_1]^*
\]
then one can identify $C^\cdot(L)$ with $K(B;f)$.

\begin{Example} Starting from $b=2$ the product
decomposition of the Hilbert series becomes inf\/inite: for example,
if $b=2$
\[
H(A_2;t) = \frac{1+2t}{(1-t)^2} = \frac
{(1-t^2)^3(1-t^4)^3(1-t^6)^{11}(1-t^8)^{30}(1-t^{10})^{105}\cdots}
{(1-t)^4(1-t^3)^2(1-t^5)^{6}(1-t^7)^{18}(1-t^{9})^{56}\cdots}.
\]
The f\/irst two exponents are ``koszul'': the number of unknowns and the
number of equations. The exponents grow exponentially.
\end{Example}

\subsection{Gauss cyclotomic identity}\label{section3.2}

\paragraph{3.2.1.~Necklaces.}
The {\it necklace polynomial} is def\/ined as
\[
M_n(x) = \frac{1}{n}\sum_{d|n} \mu(d)x^{n/d}.
\]
\begin{Example} If $p$ is a prime, then $M_p(x) = (x^p - x)/p$.
\end{Example}

Let a {\it necklace} $c$ be made of $n$ beads; suppose  that each
bead can be one of $m$ colors. A~necklace is called {\it primitive}
if it is not of the form $c = dc'$ where $d|n$, $d > 1$.


\begin{Theorem}[C.~Moreau, 1872] The number of primitive
necklaces made of $n$ beads in $b$ colors is equal to $M_n(b)$.
\end{Theorem}

The proof is an exercise. C.~Moreau was an artillery captain
from Constantine, cf.~\cite{M}.


\begin{Corollary} \label{corollary3.2.3} The number of all necklaces made of $n$
beads in $b$ colors is equal to $\Phi_n(b)$, where
\[
\Phi_n(x) = \frac{1}{n}\sum_{d|n} \phi(d)x^{n/d},
\]
$\phi(d)$ is the Euler function.
\end{Corollary}
\begin{proof} If this number is $C(n;b)$ then
\[
C(n;b) = \sum_{l|n} M_l(b) = \sum_{l|n} \sum_{d|l} \frac{\mu(d)}{d}
\cdot\frac{b^{l/d}}{l/d} \ \overset{\text{(we set $p = l/d$)}}{=} \
 \sum_{p|n} b^p\cdot \left(
\frac{1}{p}\sum_{d|(n/p)}\frac{\mu(d)}{d}\right)
\]
or,
\[
\sum_{d|(n/p)}\frac{\mu(d)}{d} = \frac{\phi(n/p)}{n/p},
\]
which proves the corollary.
\end{proof}

\paragraph{3.2.2. A theorem of P\'olya.}
Following Polyakov (cf.~\cite{Pol}), one can consider the same
numbers from the point of view of P\'olya theory~\cite{P}.

Suppose we are given two f\/inite sets $X$ and $Y$ as well as a weight
function $w: Y\rightarrow \mathbb{N}$. If $n=|X|\,$, without loss of
generality we can assume that $X = \{1,2,\ldots,n\}$. Consider the
set of all mappings $F = \{ f\mid f:X\rightarrow Y \}$. We can
def\/ine the weight of a function $f\in F\,$ to be
\[
 w(f) = \sum_{x\in X} w\left(f(x)\right).
\]

Every subgroup of the symmetric group on $n$ elements, $S_{n}$,
acts on $X\,$ through permutations. If $A$ is one such subgroup,
an equivalence relation $\sim_{A}$ on $F$ is def\/ined as $f \sim_{A}
g$ $\Longleftrightarrow$ $f = g\circ a$ for some $a\in A$.

Denote by $[f] = \{ g \in F\mid f \sim_{A} g \}$ the equivalence
class of $f\,$ with respect to this equivalence relation. $[f]$ is
also called the orbit of $f$. Since each $a \in A$ acts
bijectively on $X$, then
 \[
  w(g) = \sum_{x\in X} w\left(g(x)\right) = \sum_{x\in X} w\left(g(a\circ x))\right) = \sum_{x\in X} w\left(f(x)\right) =
  w(f).
  \]

Therefore we can safely def\/ine $w([f]) = w(f)$. In other words,
permuting the summands of a sum does not change the value of the
sum.

Let $ c_k = \left|\{ y\in Y \mid w(y)=k \}\right|$ be the number of
elements of $Y$ of weight $k$.
The generating function by weight of the source objects is $ c(t) =
\sum_k c_k\cdot t^k$.
Let $ C_k = \left|\{ [f] \mid w([f])=k \}\right|$ be the number of
orbits of weight $k$.
The generating function of the f\/illed slot conf\/igurations is $C(t) =
\sum_k C_k\cdot t^k$.

\begin{Theorem}
Given all the above definitions, P\'olya's enumeration theorem
asserts that
\[
 C(t) = Z(A)\left(c(t),c\left(t^2\right),\ldots,c(t^n)\right),
 \]
where $Z(A)$ is the cycle index $($Zyklenzeiger$)$ of $A$
\[
Z(A)(t_1, t_2, \ldots, t_n) = \frac{1}{|A|} \sum_{g\in A} t_1^{j_1(g)} t_2^{j_2(g)} \cdots
 t_n^{j_n(g)}.
 \]
\end{Theorem}

Consider the group of cyclic permutations  $G \isom \BZ/n$  as a
subgroup of the symmetric group~$S_n$.
Def\/ine, after P\'olya, the  cycle index
 polynomial of $G$, in $n$ variables, as
\[
P_G(x_1,\ldots, x_n) = \frac{1}{n} \sum_{\sigma\in G} \prod_{i=1}^n
x_i^{c_i(\sigma)},
\]
where $c_i(\sigma)$ is the number of cycles of length $i$ in
$\sigma$.

In other words, one associates to each  $\sigma \in G$ a monomial.
For example, for $n=6$ there are following permutations:
\begin{gather*}
\sigma_0 = (1)(2)(3)(4)(5)(6), \quad \text{the corresponding monomial is
$x_1^6$},\\
\sigma_1 = (123456):\quad  x_6,\\
\sigma_2 = \sigma_1^2 = (135)(246):\quad x_3^2,\\
\sigma_3 = \sigma_1^3 = (14)(25)36):\quad x_2^3,\\
\sigma_4 = \sigma_1^4 = (153)(264):\quad x_3^2,\\
\sigma_5 = \sigma_1^5 = (654321):\quad x_6.
\end{gather*}

The Zyklenzeiger is equal to $P_{\BZ/6}(x) = \frac{1}{6}(x_1^6 +
x_2^3 + 2x_3^2 + 2x_6)$.

In general,
\[
P_{\BZ/n}(x_1,\ldots,x_n) = \frac{1}{n} \sum_{d|n} \phi(d) x_d^{n/d}.
\]
After a change of variables $x_i = \sum\limits_{j=1}^b y_j^i$, one obtains
a polynomial
\[
W_G(y_1,\ldots,y_b) := P_G\left(\sum y_j, \sum y^2_j, \ldots, \sum y_j^n\right).
\]
As a special case of the theorem of P\'olya we obtain that the
number of necklaces made of $n$ beads in $b$ colors is equal to
$W_{\BZ/n}(1,\ldots,1) = P_{\BZ/n}(b,\ldots,b) = \Phi_n(b)$.

\paragraph{3.2.3.\ Cyclotomic identity.}
\begin{Theorem}[Gauss, cf.~\cite{G}] One has the following formula:
\[
1 - bt = \prod_{n=1}^\infty  \left(1 - t^n\right)^{M_n(b)}.
\]
\end{Theorem}
It is proved by the application of the M\"obius inversion. There is a
useful generalization of this identity, found by Pieter Moree~\cite{Mor}. We will need it later studying the equivariant Hilbert
series, see Section~3.4.4.


\begin{Theorem} \label{theorem3.2.6}
Let $f(q) \in \BC[q,q^{-1}]$ be an
arbitrary Laurent polynomial. Introduce the polynomials
\[
M_n(f;q) = \frac{1}{n}\sum_{d|n} \mu(d)f(q^d)^{n/d}.
\]
Then
\[
1 - f(q)t = \prod_{n=1}^\infty \prod_{i=-\infty}^\infty \left(1 -
q^it^n\right)^{a_{in}},
\]
where for each $n$, the number $a_{in}$ are defined by the equations
\[
\sum_i a_{in}q^i = M_n(f;q).
\]
\end{Theorem}

\begin{proof} Applying $-\log$ to the both sides of the above
formula we obtain:
\[
-\log(1 - f(q)t) = \sum_{m=1}^\infty  \frac{f(q)^mt^m}{m}
\]
and
\begin{gather*}
-\log\left(\prod_{i,n}\big(1 - q^it^n\big)^{a_{in}}\right) = \sum_i\sum_{n,k=1}^\infty
a_{in}\frac{q^{ik}t^{nk}}{k} =
 \sum_{m=1}^\infty t^m\cdot\left(\sum_{n|m}\sum_i a_{in}
\frac{q^{im/n}}{m/n}\right),
\end{gather*}
therefore
\[
f(q)^m = \sum_{n|m}\sum_i na_{in}q^{im/n}
\]
for each $m=1,2,\ldots $. Make a change of variables: $p = q^m$,
\[
f(p^{1/m})^m = \sum_{n|m}\sum_i na_{in}p^{i/n}.
\]
Using the M\"obius inversion we get:
\[
\sum_i na_{in}p^{i/n} = \sum_{d|n} \mu(n/d) f(p^{1/d})^d,
\]
and making a substitution: $q=p^{1/n}$, we obtain
\[
\sum_i na_{in}q^{i} = \sum_{d|n} \mu(n/d) f(q^{n/d})^d,
\]
which is the formula required.
\end{proof}

Moreover, O.~Ogievetsky noticed that the theorem can be generalized
to the case of several variables:

\begin{Theorem}\label{theorem3.2.7}
 Let $f(q_1,\ldots,q_p)$ be a Laurent
polynomial in $p$ variables. Introduce the  Laurent polynomials
\[
M_n(f;q_1,\ldots,q_p) = \frac{1}{n}\sum_{d|n}
\mu(d)f\big(q_1^d,\ldots,q_p^d\big)^{n/d}.
\]
Then
\[
1 - f(q,\ldots,q_p)t = \prod_{n=1}^\infty \prod_{i=-\infty}^\infty
\big(1 - q_1^{i_1}\cdots q_p^{i_p} t^n\big)^{ a_{i_1\ldots i_p;n}},
\]
where for each $n$, the numbers $a_{i_1\ldots i_p;n}$ are the
coefficients of $M_n$:
\[
\sum_i  a_{i_1\ldots i_p;n}q_1^{i_1}\cdots q_p^{i_p} =
M_n(f;q_1,\ldots, q_p).
\]
\end{Theorem}
The proof is the same as before.


\begin{Example} \label{example3.2.8} Take $f(q) = - q$. Then
\[
1 + qt = (1 - qt)^{-1}\big(1 - q^2t^2\big),
\]
where
\[
\frac{1}{n}\sum_{d|n} (-1)^d\mu(d)
\]
is equal to $0$ if $n \geq 3$, $-1$ if $n=1$, and $1$ if $n=2$.
\end{Example}

\paragraph{3.2.4. Application: the $\boldsymbol{\zeta}$ function of the af\/f\/ine line.}
Let $A := \BF_p[x]$; this ring is similar in many ways to $\BZ$.

Nonzero ideals $I \subset A$ are in bijection with unitary
polynomials $f(x)$, $I = (f)$, and principal ideals correspond to irreducible
ones. Set
\[
N(I) := \sharp (A/I) = p^{\deg f},
\]
and def\/ine
\[
\zeta(A;s) = \sum_{I\subset A,\, I\neq 0} N(I)^{-s} = \sum_{f\
\text{unitary}} p^{- s\deg f}.
\]
There is $p^n$ unitary polynomials of degree $n$, therefore
\begin{gather}
\zeta(A;s) = \sum_{n=1}^\infty p^n \cdot p^{-sn} = \frac{1}{1 -
p\cdot p^{-s}} = \frac{1}{1 - pT}, \label{(3.2.9.1)}
\end{gather}
where $T := p^{-s}$.

The Euler product formula for $\zeta(A;s)$ can be written in the
following form
\begin{gather*}
\zeta(A;s) = \prod_{f\ \text{unitary, irreducible}} \frac{1} {1 -
p^{-\deg f\cdot s}}  \\
\phantom{\zeta(A;s)}{}
= \prod_{d=1}^\infty \prod_{f\ \text{un., irr.}, \, \deg f= d}
\frac{1}{1 - p^{-ds}} = \prod_{d=1}^\infty
 \frac{1}{(1 - T^d)^{N_d(p)}},
\end{gather*}
where $N_d(p)$  denotes the number of unitary irreducible
polynomials of degree $d$ in $A$.

On the other hand, applying the cyclotomic identity to \eqref{(3.2.9.1)}, one
gets
\[
\zeta(A;s) = \frac{1}{1 - pT} = \frac{1}{\prod\limits_{d=1}^\infty (1 -
T^d)^{M_d(p)}},
\]
proving therefore

\begin{Theorem}[Gauss~\cite{G}] The number $N_d(p)$ of irreducible
unitary polynomials of degree $d$ in~$\BF_p[x]$ is equal to
\[
M_d(p) = \frac{1}{d}\sum_{l|d}\mu(l)p^{d/l}.
\]
\end{Theorem}
\begin{Corollary} For $d\geq 1$, $N_d(p) > 0$, i.e.\ for each $d\geq
1$ there is an irreducible polynomial of degree~$d$.
\end{Corollary}
It is interesting to compare this theorem of Gauss with Riemann's
explicit formula:
\[
\pi(x) = \sum_{n=1}^\infty
\frac{\mu(n)}{n}\left\{\text{Li}(x^{1/n}) + \sum_{{\rm {\frak I}m}(\rho)>0}  (x^{\rho/n} + x^{1-\rho}) + \int_{x^{1/n}}^\infty
\frac{dt}{t(t^2-1)\log t} - \log 2\right\}.
\]
Here $\pi(x) = $ the number of primes $p \leq x$, the summation is
over the non-trivial roots of  $\zeta(s)$, and
\[
\text{Li}\, x = \int_1^x  \frac{dt}{\log t}.
\]

\subsection{Witt formula}\label{section3.3}

\begin{flushright}
\begin{minipage}{8cm}
\it La Nature est un temple o\`u de vivants piliers\\
Laissent parfois sortir de confuses paroles \dots

\vspace{2mm}
\quad  \hfill Ch.~Baudelaire
\end{minipage}
\end{flushright}

\smallskip

\begin{Theorem}[Witt~\cite{W}] \label{theorem3.3.1} Let $L$ be a free Lie algebra
on $b$ generators. Then
\[
\dim L_n = M_n(b) = \frac{1}{n}\sum_{d|n} \mu(d)b^{n/d}.
\]
In other words,
\[
H(L;t) := \sum_{n=0}^\infty \dim L_n\cdot t^n = - \sum_{m=1}^\infty
\frac{\mu(m)}{m}\log \big(1 - bt^m\big).
\]
\end{Theorem}
\begin{proof} If $x_1,\ldots,x_b$ are the generators of $L$,
then the universal enveloping algebra $UL = \BC\langle x_1,\ldots,x_b\rangle$ (a
free associative algebra), has the following Hilbert series
\[
H(UL;t) := \sum_{i=0}^\infty \dim UL_i\cdot t^i = \frac{1}{1 - bt}.
\]
On the other hand, due to Poincar\'e--Birkhof\/f--Witt, $UL = S^\cdot L$,
therefore, if one denotes $a_n := \dim L_n$, then
\[
H(UL;t) = \prod_{n=1}^\infty \frac{1}{(1 - t^n)^{a_n}},
\]
and applying the cyclotomic identity f\/inishes the proof.
\end{proof}

\paragraph{3.3.1.}
This theorem admits a nice generalization. Let $L$ be a
free Lie superalgebra with $b$ even generators $x_1,\ldots, x_b$ and
$c$ odd generators  $y_1,\ldots, y_c$.

Equip $L$ with a grading, assigning to $x_i$ and to $y_j$ the
degree $1$, and setting $\deg([a,b]) = \deg(a) + \deg(b)$. If $L_n
\subset L$ is a subspace of the elements of the degree $n$, then $L
= \mathop{\oplus}\limits_{n=1}^\infty L_n$.

Explicitly, the Lie monomial  $\alpha = [a_1,[a_2,\ldots
a_n]\ldots]$, where $a_i = x_j$ or $y_k$, has the degree $n$. On the
other hand, its {\it parity} is equal to the number of $y$'s among
$a_i$. For example, $[x_i,y_j]$ is an element odd of degree $2$
(sic!).

Each homogenous summand $L_n$ therefore is decomposed into two
summands: even and odd elements, $L_n = L_n^p \oplus L_n^i$. One
def\/ines the {\it super-dimension} as $\dim^\sim L_n = \dim L_n^p -
\dim L_n^i$.

\begin{Theorem}[Petrogradsky~\cite{Pet}] One has the formula
\begin{gather}
\dim^\sim L_n = M_n(b-c) = \frac{1}{n}\sum_{d|n} \mu(d)(b-c)^{n/d}.
\label{(3.3.2.1)}
\end{gather}
\end{Theorem}
\begin{Example} Let $b=0$, $c=1$; then $L = \BC\cdot y \oplus \BC\cdot
[y,y]$. In the right hand side, one sees that $M_n(-1) = -1$ if
$n=1$, $1$ if $n=2$, and $0$ if $n\geq 3$, cf.\ Examples~\ref{example3.2.8} and~\ref{example3.1.4}.
\end{Example}

\paragraph{3.3.2.\ Hochschild homology.}
Let $A$ be an  associative unitary algebra. Recall that the
Hochschild homology  $HH_i(A)$  of $A$ is def\/ined as the homology of
the complex
\[
CH_\cdot (A):\ \cdots \lra A\otimes A\otimes A \buildrel b\over\lra
A\otimes A \buildrel b\over\lra A  \lra 0,
\]
where $CH_n(A) = A^{\otimes n+1}$, and the dif\/ferential $b$ is
def\/ined by
\begin{gather*}
b(a_1\otimes a_2\cdots\otimes a_n) = a_1a_2\otimes a_3\otimes \cdots
\otimes a_{n} - a_1\otimes a_2a_3\otimes a_4 \cdots\otimes a_n +\cdots
\\
\qquad{} + (-1)^{n-1}a_1\otimes\cdots\otimes a_{n-1}a_n + (-1)^n
a_na_1\otimes a_2\otimes \cdots\otimes a_{n-1}.
\end{gather*}
In particular, $d(a\otimes b) = ab - ba$ and $HH_0(A) = A/[A,A]$,
where $[A,A] \subset A$ is the subspace generated by the elements
$ab - ba$.

In fact, one can replace $A$ by $\bar A := A/\BC$: the complex
$CH_\cdot (A)$ is quasi-isomorphic to
\[
\bar{CH}_\cdot (A):\ \cdots \lra A\otimes \bar A\otimes \bar A \lra \bar
A\otimes A \lra A  \lra 0.
\]
One has
\[
HH_i(A) = {\rm Tor}^{A\otimes A^o}_i(A,A).
\]
Let us introduce the operators $B:  A^{\otimes n+1} \lra A^{\otimes
n+2}$ by the formula:
\begin{gather*}
B(a_0,\ldots,a_n) = \sum_{i=0}^{n-1}
\big\{(-1)^{ni}(1,a_i,\ldots,a_n,a_0,\ldots, a_{i-1})\\
\phantom{B(a_0,\ldots,a_n) =}{}  -
(-1)^{n(i-1)}(a_{i-1},1,a_i,\ldots,a_n,a_0,\ldots, a_{i-2})\big\}.
\end{gather*}
For example,
\begin{gather}
B(a_0) = (1,a_0) + (a_0,1). \label{(3.3.3.1)}
\end{gather}
Then $Bb + bB = 0$ (cf.~\cite[2.1]{L}); therefore it induces the
morphisms
\begin{gather}
B:\ HH_n(A) \lra HH_{n+1}(A). \label{(3.3.3.2)}
\end{gather}
Recall that the cyclic homology can be def\/ined as the homology of
the bi-complex
\begin{gather*}
HC_i(A) = H_i\bigl(CH_\cdot (A)\buildrel B\over\lra CH_\cdot (A)[1] \buildrel
B\over\lra CH_\cdot (A)[2] \buildrel B\over\lra \cdots\bigr).
\end{gather*}


\paragraph{3.3.3.} Now one can give an algebraic interpretation of the
polynomials $\Phi_n(x) = \sum\limits_{d|n}\phi(d)x^{n/d}/n$.

Let $A = \BC\langle x_1,\ldots,x_b\rangle$, then $HH_0(A)$ inherits  a grading
from $A$, $HH_0(A) = \mathop{\oplus}\limits_{n=0}^\infty HH_0(A)_n$; if one has
$\bar A = \mathop{\oplus}\limits_{n\geq 1}A_n$, then:
\[
\bar{HH}_0(A) = \mathop{\oplus}\limits_{n\geq 1} HH_0(A)_n = \bar A/[\bar A,\bar A].
\]

One can think of $HH_0(A)_n$ as a space of {\it cyclic words}
$x_{i_1}x_{i_2}\ldots x_{i_n}$ of length $n$ in letters $x_i$, where
two words are identif\/ied if one is a cyclic permutation of another.
Therefore cyclic words are identif\/ied with necklaces of $n$ beads in
$b$ colors.

Hence $HH_0(A)_n$ can be viewed as a linear space with a basis
indexed by necklaces of $n$ beads in $b$ colors. It follows that
\[
\dim HH_0(A)_n = \Phi_n(b) = \frac{1}{n}\sum_{d|n}\phi(d)b^{n/d},\qquad
n\geq 1,
\]
(cf.\ Corollary~\ref{corollary3.2.3}), or
\[
H(\bar{HH}_0(A);t) = - \sum_{m=1}^\infty \frac{\phi(m)}{m}\log(1 -
bt^m)
\]
({\it the Polyakov  formula}).

\paragraph{3.3.4.}
Let $V = \mathop{\oplus}\limits_{i=1}^b \BC\cdot x_i$, therefore $A = TV$.
For each $n\geq 1$ one def\/ines an automorphism $\tau:  V^{\otimes n}
\lra V^{\otimes n}$ as
\[
\tau(v_1\otimes\cdots\otimes v_n) = v_n\otimes v_1\otimes
v_2\otimes\cdots \otimes v_{n-1}.
\]
One observes that
\[
HH_0(TV)_n = V^{\otimes n}_{\tau} := \text{Coker}(1-\tau).
\]
Following \cite[3.1]{L}, one def\/ines a complex of length $1$:
\[
CH_\cdot ^{sm}(TV):\  0 \lra TV\otimes V \lra TV \lra 0
\]
equipped with a dif\/ferential $d(a\otimes v) = av - va$. Def\/ine also
a morphism of complexes $\phi: CH_\cdot (TV) \lra CH_\cdot ^{sm}(TV)$ as:
$\phi_0 = \text{Id}_{TV}$,
\begin{gather}
\phi_1(a\otimes v_1\ldots v_n) = \sum_{i=1}^n v_{i+1}\cdots v_n a
v_1\cdots v_{i-1}\otimes v_i. \label{(3.3.5.1)}
\end{gather}
On the other hand, there is an evident inclusion  $\iota:
CH_\cdot ^{sm}(TV) \lra CH_\cdot (A)$ such that $\phi\circ\iota = \text{Id}$,
and, as one can check, $\iota\circ\phi$ is homotopic to the
identity.

It follows that
\[
HH_1(TV)_n = \big(V^{\otimes n}\big)^{\tau} := \text{Ker}(1-\tau)
\]
and $HH_i(TV) = 0$ for $i\geq 2$.

Instead of the Hochschild homology one can also consider the
(reduced) cyclic homology. Then: $HC_0(TV) = HH_0(TV)$ and
$\bar{HC}_i(TV) = 0$ for $i>0$, cf.~\cite{Mov}.

\paragraph{3.3.5.~Partial derivatives.}
Let $m = \dots x_ix_jx_k\dots$ be a
cyclic word (i.e.~$m\in TV/[TV,TV]$) such that the letter $x_i$
appears once in it; then one can def\/ine a usual word (i.e.\ an
element de~$TV$) $\dpar m/\dpar x_i$, by ``cutting'' $m$ and
removing the letter $x_i$:
\[
\dpar m/\dpar x_i = x_jx_k \dots\,.
\]
When $x_i$ appears in a word several times, the result will be a
sum, by the Leibnitz rule, cf.~\cite{K}. In this way we def\/ine a map
\[
\frac{\dpar}{\dpar x_i}:\ TV/[TV,TV] \lra  TV.
\]
Consider an operator $\sum\limits_{i=1}^b x_i\dpar/\dpar x_i$ (cf.~\eqref{(3.3.5.1)}). It respects the polynomial degree, and it is not hard to
verify that its image is contained in $TV^\tau$. For example:
\[
e(xyz) = xyz + yzx + zxy.
\]
One observes that the map
\begin{gather}
\sum_{i=1}^b x_i\frac{\dpar}{\dpar x_i}:\ HH_0(TV) = TV_\tau \lra
TV^\tau = HH_1(TV) \label{(3.3.6.1)}
\end{gather}
coincides with homomorphism $B$,  cf.~\eqref{(3.3.3.1)}, \eqref{(3.3.3.2)}. (Compare
\cite[the line before (14), p.~8]{Mov}.)

It was pointed out to us by V.~Ginzburg that there is another
interpretation of the above map. Consider the composition
\begin{gather}
d:\ TV_\tau = TV/[TV,TV] \buildrel{\sum x_i\dpar/\dpar x_i}\over\lra
TV^\tau \subset TV \otimes V . \label{(3.3.6.2)}
\end{gather}
It was already Quillen who studied this map in the $80$'s, cf.~\cite{Q2}. In the notations of \cite[\S~3]{Q2} $TV/[TV,TV] =
TV_\natural$ and $TV \otimes V = \Omega^1_{TV,\natural}$; one can
consider these spaces as the space of cyclic functions
(dif\/ferentiable $1$-forms respectively)  over the ``non-commutative
space $\text{Spec}\ TV$''; the morphism \eqref{(3.3.6.2)} is called ``the Karoubi--de Rham dif\/ferential''.

\subsection{Koszul duality}\label{section3.4}

\paragraph{3.4.1.} Let $A$ be an associative {\it quadratic} algebra, which
means that, it is a quotient $A = TV/(R)$ of the free associative
algebra over the space $V$ of f\/inite dimension by the two-sided
ideal
 $(R)$ generated by a subspace of relations $R
\subset V\otimes V$.

Recall that the quadratic dual algebra is def\/ined as $A^! =
TV^*/(R^\perp)$, where $R^\perp \subset V^*\otimes V^* \isom
(V\otimes V)^*$ is the annihilator of $R$.

We are interested in the case when  $A$ is commutative; in this case
$R \supset \Lambda^2V$, and $R^\perp \subset \Lambda^2V^\perp =
S^2V^*$. In other words, if $\{x_i\}$ form a basis of $V$, $\{y_i\}$
form a basis of $V^*$, then $R^\perp$ is contained in $S^2V^*
\subset V^{\otimes 2}$; in other words it is generated by odd
commutators $[y_i,y_j] = y_iy_j + y_iy_j$. Def\/ine a Lie algebra $L$
as a {\it graded} Lie algebra generated by $y_i$ of degree $1$ and
relations $g=0$, $g\in R^\perp$. Then $A^! = UL$ by def\/inition. The
Lie algebra $L$ is called the {\it Koszul, or Quillen dual} of $A$.

For example, if $A = H^*(X;\BC)$ is the cohomology ring of a simply
connected topological space~$X$, then $L$ is its homotopy Lie
algebra, $\oplus \, \pi_i(X)_\BC$, under some additional conditions,
cf.~\cite{Q1}.

Therefore we have, $L = \mathop{\oplus}\limits_{n=1}^\infty L_n$; moreover one
observes that $L_1 = V^*$, $L_2 = S^2V^*/R^\perp$. Set $L_{\leq 2} =
L_1 \oplus L_2$; it is a quotient Lie algebra of $L$.

\paragraph{3.4.2.} {\it The Chevalley cochain complex} of $L$ is by def\/inition
the space
\[
C^\cdot(L) = SL^*_{1}\otimes \Lambda L^*_{2}\otimes
SL^*_{3}\otimes\cdots
\]
equipped with the Chevalley dif\/ferential. This complex is double
graded:
\[
C^\cdot(L) = \mathop{\oplus}\limits_{a,b=0}^\infty C^\cdot(L)_{ab},
\]
where the f\/irst (homological) degree of $L_i^*$ is set to be equal
$i-1$, and the second (polynomial) degree of $L_i^*$ is $i$, both
gradings are compatible with the product. The Chevalley
dif\/ferential preserves the second grading and decreases the f\/irst by~$1$:
\[
dC^\cdot(L)_{ab}
 \subset C^\cdot(L)_{a-1,b}.
\]
Here are the components of $C^\cdot(L)$  of polynomial
degree $\leq 3$:
\begin{gather}
L^*_3 \lra L_2^*\otimes L_1^* \lra S^3L_1^*,\nonumber
\\
L_2^* \lra S^2L_1^*, \label{(3.4.2.1)}
\\
L_1^*.\nonumber
\end{gather}
One has: $C^\cdot(L)_{0\cdot} = SL_1^*$ and the complex starts as:
\begin{gather}
\cdots\lra (L_3^*\oplus \Lambda^2L_2^*)\otimes SL_1^* \lra
L_2^*\otimes SL_1^* \lra SL_1^* . \label{(3.4.2.2)}
\end{gather}
Since $L_1^*= V$ and $L_2^* = R\subset S^2V$, the f\/irst dif\/ferential
in~\eqref{(3.4.2.2)} is a $SV$-linear map which sends $f\otimes 1\in L_2^*$
to $f\in R\subset S^2L_1^*$.

One observes that:

(a) there is a natural augmentation $C^\cdot(L) \lra A = SV/(R)$,
cf.~\cite{MS}.

(A similar picture arises in the construction of the mixed Tate
motivic cohomology of a~f\/ield~$k$, cf.~\cite{BGSV}. There, the
analogues of complexes $C^\cdot(L)_{n\cdot}$ are the Beilinson
motivic complexes $\BZ(n)$, the analogue of $A$ is the Milnor
$K$-theory $K_\cdot ^{\rm Miln}(k)$; the analogue of  $L$ is (the Lie algebra
of) the ``mixted Tate motivic fundamental group''.)

(b) If one chooses a basis $\{f_1,\ldots,f_d\}$ of $R$, then one can
identify the Chevalley complex of  $L_{\leq 2}$ with the Koszul
complex
\[
C^\cdot (L_{\leq 2}) \isom K(SV;(f)),
\]
cf.\ Example~\ref{example3.1.4}. 

(c) There is a natural inclusion $C^\cdot (L_{\leq 2}) \subset
C^\cdot (L)$ compatible with the augmentation to $A$.

Def\/ine a generating series
\[
H(C^\cdot(L);u,t) := \sum_{a,b = 0}^\infty  \dim  C^\cdot(L)_{ab}
u^at^b
\]
and its specialization, the Euler--Hilbert series:
\[
EH(C^\cdot(L);t) := H(C^\cdot(L);-1,t) = \sum_{n=0}^\infty
EP(C^\cdot(L)_{?,n}) t^n,
\]
where $EP$ stands for the Euler--Poincar\'e characteristic. One
observes immediately that
\[
EH(SL^*_n) = (1 - t^n)^{-\dim L_n}\qquad \text{and} \qquad EH(\Lambda L^*_n) =
(1 - t^n)^{\dim L_n},
\]
where
\begin{gather}
EH(C^\cdot(L);t) = \prod_{n=1}^\infty (1 - t^n)^{(-1)^n\dim L_n}.
\label{(3.4.2.3)}
\end{gather}

\paragraph{3.4.3.} For example, consider the Veronese ring $A_b$, cf.\ Section~3.1.1.
Recall that $A_b$ is generated by $x_0, \ldots, x_{b+1}$, subject to
relations $x_ix_j - x_kx_l = 0$ if $i+j=k+l$. Then $A^!$ is
generated by $y_i = x_i^*, 0 \leq i \leq b+1$ obeying the relations
\begin{gather}
\sum_{i+j=k} [y_i,y_j] = 0,\qquad k=0,\ldots,2(b+1), \label{(3.4.3.1)}
\end{gather}
where $[y_i,y_j] = y_iy_j + y_jy_i$ (sic!). Therefore if one def\/ines
a graded Lie algebra $L^b$ with generators $y_i \in L^{b}_1 =
A_{b1}^*$ and relations~\eqref{(3.4.3.1)}, then $A_b^! = UL^b$.

This Lie algebra has a nice structure, cf.~\cite{GKR} (the case
$b=1$ was considered in Example~\ref{example3.1.4}). It admits an involution $i(y_j) =
y_{b+1-j}$. Let ${\hat L}^b \subset L^b$ be a Lie subalgebra
generated by $y_1,\ldots,y_b$. Then ${\hat L}^b = [L^b,L^b]$ is a
Lie ideal (stable under $i$); as a graded Lie algebra it is {\it
free}.

The quotient algebra $\bar L = L^b/{\hat L}^b$ is a graded Abelian
Lie algebra on $2$ generators~$\bar y_0$,~$\bar y_{b+1}$. One sees
that ${\hat L}^b_{\geq 2} = L^b_{\geq 2}$.

We have noticed in Example~\ref{example3.1.4} that $L_1$ is f\/inite-dimensional; as
opposed to that, if $b \geq 2$, than $L^b$ is inf\/inite-dimensional
(with an exponential growth). Namely, one has $L^b =
\mathop{\oplus}\limits_{n=1}^\infty L_{bn}$, where the dimensions of the homogenous
components can be calculated easily using an odd analogue of the
Witt theorem, see Theorem~\ref{theorem3.3.1}. One obtains:
\begin{gather*}
\dim L_1 = - M_1(-b) + 2 = b+2,
\\
\dim^\sim L_n = \dim^\sim \hat L_n = M_n(-b) = \frac{1}{n}
\sum_{d|n} \mu(d) (-b)^{n/d}\qquad (n\geq 2).
\end{gather*}
Here $\hat L_{b1}$ is generated by $b$ odd elements, therefore
\[
\dim^\sim \hat L_n = (-1)^n \dim \hat L_n.
\]
One observes that the signs od $\dim^\sim L_n$ alternate. For
example, $\dim L_2 = M_2(-b) = (b^2 + b)/2$,  $\dim L_3 = - M_3(-b) =
(b^3 - b)/3$, etc.

Now consider the Chevalley complex of $L^b$: \eqref{(3.4.2.3)} implies that
its Euler--Hilbert is as follows:
\begin{gather}
EH\big(C^\cdot\big(L^b\big);t\big) = (1-t)^{-\dim L_{b1}}\left(1-t^2\right)^{\dim L_{b2}}
\left(1-t^3\right)^{-\dim L_{b3}}\cdots  \nonumber\\
\phantom{EH\left(C^\cdot\left(L^b\right);t\right)}{}
= (1-t)^{-2}\prod_{n=1}^\infty\left(1-t^n\right)^{M_n(-b)} = \frac{1+bt}{(1 -
t)^2}. \label{(3.4.3.2)}
\end{gather}
It is clear that it coincides with the Hilbert series of the ring
$A_b$, \eqref{(3.1.3.1)}. the last equality is due to the cyclotomic
identity. It is not surprising.

In fact, a deep theorem of Bezrukavnikov says:

\begin{Theorem}[cf.~\cite{Bez}] \label{theorem3.4.4} The algebra $A_b$ is Koszul.
\end{Theorem}

Formula~\eqref{(3.4.3.2)} can be viewed as a  ``numerical evidence'' that
the theorem holds.

\begin{Corollary} There are natural isomorphisms $A_b^! \isom
{\rm Ext}^\cdot _{A_b}(\BC,\BC)$, $A_b\isom {\rm Ext}^\cdot _{A_b^!}(\BC,\BC) =
H^\cdot (L;\BC)$.
\end{Corollary}
It follows that $A_b^\sim := C^\cdot (L^b)$ is a dga resolution of $A_b$,
cf.~\cite{MS}.

\paragraph{3.4.4.~Characters.}
The Lie algebra $\text{sl}(2)$ acts on
$A_b$ in such a way that $A_{b1}$ is an irreducible
$\text{sl}(2)$-module; therefore its character is
\[
{\rm Ch}(A_{b1}) = [b+1]_q = \frac{q^{b+1} - q^{-b-1}}{q - q^{-1}}.
\]
The character of $A_b$ is given by the equivariant Hilbert series,
cf.~\eqref{(3.1.3.1)_q}:
\[
H_q(A_b;t) = \frac{1 + [b]_qt}{(1 - q^{b+1}t)(1 - q^{-b-1}t)}.
\]
The above action induces an action of $\text{sl}(2)$ on $L^b$ and
therefore on $A_b^\sim := C^\cdot (L^b)$. Applying now the
``$q$-cyclotomic'' identity: Theorem~\ref{theorem3.2.6} for $f(q) = - [b]_q$,
one obtains the $\text{sl}(2)$-character ${\rm Ch}(A_b^\sim)$.

On the other hand:the subalgebra Lie $\hat L_b \subset L^b$ is free
on $b$ generators, therefore the Lie algebra $\text{gl}(b)$ acts on
it. One can use the theorem of Ogievetsky (see Theorem~\ref{theorem3.2.7}) for
$f(q_1,\ldots,q_b) = - {\rm Ch}_{\text{gl}(b)}(\hat L_{b1})$ for
calculation of the character ${\rm Ch}_{\text{gl}(b)}(C^\cdot (\hat L_b))$.

The $\text{gl}(b)$-character of the free Lie algebra on {\it odd}
$b$ generators was calculated by Angeline Brandt in \cite{Br}.

The free {\it group} on $b$ generators is isomorphic to the
fundamental group of the Riemann sphere with $b+1$ points removed;
its nilpotent completion is a fundamental object of the theory of
Grothendieck--Drinfeld--Ihara, cf.~\cite{D}.

\paragraph{3.4.5.}
Returning to the case of an arbitrary commutative quadratic
algebra, one can show that a complete intersection is Koszul, cf.~\cite{PP}. The other way around, if $A$ is commutative Koszul, and
$L$ is its dual Lie algebra, then $A$ is a complete intersection if
and only if $L = L_{\leq 2}$.

One can say that commutative (or maybe also noncommutative?) Koszul algebras
are natural generalizations of the quadratic complete intersections; and
one would expect that all the results which hold for quadratic
complete intersections will generalize to Koszul algebras.

\subsection*{Acknowledgements}
We thank Fedor Malikov who read thoroughly the f\/irst part and helped
to correct many signs; some calculations made with him have been the
starting point of the second part. We are grateful to Vladimir
Hinich for very interesting discussions about Golod rings and Koszul
algebras; to Oleg Ogievetsky for an important remark; to Alexander
Polishchuk for very useful consultations; to Hossein Abbaspour and
Thomas Tradler who taught us about the string topology, and
especially to Victor Ginzburg for his numerous explanations,
questions and bibliographical comments.

This article was f\/inished during our stay at Max-Planck-Institut
f\"ur Mathematik and the Hausdorf\/f Institut f\"ur Mathematik in June
and July 2008; we are grateful to both institutions for the
excellent working atmosphere.

\pdfbookmark[1]{References}{ref}
\LastPageEnding


\begin{thebibliography}{99}

\footnotesize\itemsep=0pt


\bibitem{AABN} Aisaka Yu., Arroyo E.A., Berkovits N., Nekrasov N.,
 Pure spinor partition function and the massive superstring spectrum, {\it J. High Energy Phys.} {\bf 2008}  (2008),  no.~8, 050, 72~pages, \href{http://arxiv.org/abs/0806.0584}{arXiv:0806.0584}.

\bibitem{BL} Beauville A., Laszlo Y., Conformal blocks and generalized theta functions,
{\it Comm. Math. Phys.}  {\bf 164} (1994), 385--419, \href{http://arxiv.org/abs/alg-geom/9309003}{alg-geom/9309003}.

\bibitem{BD} Beilinson A., Drinfeld V., Chiral algebras, {\it American Mathematical Society Colloquium Publications}, Vol.~51, American Mathematical Society, Providence, RI, 2004.

\bibitem{BGSV}Beilinson A.A., Goncharov A.B., Schechtman V.V., Varchenko A.N.,
Aomoto dilogarithmes, mixed Hodge structures, and motivic cohomology
of pairs of triangles on the plane,
in The Grothendieck Festschrift, Vol.~I, {\it Progr. Math.}, Vol.~86, Birkh\"auser Boston, Boston, MA, 1990, 135--172.


\bibitem{BN} Berkovits  N., Nekrasov  N.A., The character of pure
spinors,  {\it Lett. Math. Phys.} {\bf  74}  (2005),  75--109, \mbox{\href{http://arxiv.org/abs/hep-th/0503075}{hep-th/0503075}}.

\bibitem{Bez}Bezrukavnikov R., Koszul property and Frobenius splitting of Schubert
varieties, \href{http://arxiv.org/abs/alg-geom/9502021}{alg-geom/9502021}.

\bibitem{CJ} Cohen R.L., Jones J.D.S., A homotopy theoretic realisation
of string topology, {\it Math. Ann.} {\bf 324} (2002), 773--798, \href{http://arxiv.org/abs/math.GT/0107187}{math.GT/0107187}.

\bibitem{CS}Chas M., Sullivan D., String topology, \href{http://arxiv.org/abs/math.GT/9911159}{math.GT/9911159}.


\bibitem{Br}  Brandt A.J., The free Lie ring and Lie representations
of the full linear group, {\it Trans. Amer. Math. Soc.} {\bf  56} (1944), 528--536.

\bibitem{D}Drinfeld V.G., On quasitriangular quasi-Hopf algebras and on a group that is closely connected with $\mathrm{Gal}(\overline{\mathbb Q}/\mathbb Q)$,
{\it Algebra i Analiz} {\bf 2} (1990), no.~4, 149--181 (English transl.: {\it Leningrad Math. J.} {\bf 2} (1991), no.~4, 829--860).

\bibitem{G}Gauss C.F., Disquisitiones generales de congruentis, Analysis residuorum.
Caput octavum, Collected Works, Vol.~2, Georg Olms Verlag,  Hildersheim~-- New York,  1973, 212--242.

\bibitem{Gol}Golyshev V., The canonical strip.~I, \href{http://arxiv.org/abs/0903.2076}{arXiv:0903.2076}.

\bibitem{GKR}Gorodentsev A.L., Khoroshkin A.S., Rudakov A.N.,  On syzygies of
highest weight orbits, {\it Amer. Math. Soc. Transl. Ser. 2}, Vol.~221, Amer. Math. Soc., Providence, RI, 2007,
79--120, \href{http://arxiv.org/abs/math.AG/0602316}{math.AG/0602316}.

\bibitem{Ha} Hardy G.H., Divergent series, The Clarendon Press, Oxford, 1949.

\bibitem{H1} Hirzebruch F., Neue topologische Methoden in der algebraischen
Geometrie, {\it Ergebnisse der Mathematik und ihrer Grenzgebiete (N.F.)}, Heft~9, Springer-Verlag, Berlin~-- G\"ottingen~-- Heidelberg, 1956.

\bibitem{H2} Hirzebruch F., Private communication, 2008.

\bibitem{L}Loday J.-L., Cyclic homology, Springer-Verlag, Berlin, 1992.

\bibitem{Ka}  Kaufmann R.M., A proof of a cyclic version of Deligne's conjecture
via cacti, \href{http://arxiv.org/abs/math.QA/0403340}{math.QA/0403340}.

\bibitem{K} Kontsevich M., Formal (non)commutative symplectic geometry,
{\it The Gelfand Mathematical Seminars, 1990--1992}, Editors L.~Corwin et al.,
Birkh\"auser Boston, Boston, MA, 1993, 173--188.

\bibitem{LZ} Lian B.H., Zuckerman G.J., New perspectives on the BRST-algebraic structure
of string theory, {\it Comm. Math. Phys.} {\bf 154} (1993), 613--646, \href{http://arxiv.org/abs/hep-th/9211072}{hep-th/9211072}.


\bibitem{M}Moreau C., Sur les permutations circulaires distinctes, {\it Nouv. Ann.
Math.}  {\bf 11} (1872), 309--314.

\bibitem{Mor} Moree P., On the average number of elements in a f\/inite f\/ield with
order or index in a prescribed residue class, {\it Finite Fields Appl.} {\bf 10} (2004), 438--463, \href{http://arxiv.org/abs/math.NT/0212220}{math.NT/0212220}.

\bibitem{Mov}Movshev M., On deformations of Yang--Mills algebras, \href{http://arxiv.org/abs/hep-th/0509119}{hep-th/0509119}.

\bibitem{MS}Movshev M., Schwarz A., Algebraic structure of Yang--Mills theory,
\href{http://arxiv.org/abs/hep-th/0404183}{hep-th/0404183}.

\bibitem{Pet}Petrogradsky V.M., On Witt's formula and invariants for free Lie
superalgebras, in Formal Power Series and Algebraic
Combinatorics (Moscow, 2000), Springer, Berlin,  2000, 543--551.

\bibitem{PP}Polishchuk A., Positselski L., Quadratic algebras, {\it University Lecture Series}, Vol.~37, American Mathematical Society, Providence, RI, 2005.

\bibitem{P}P\'olya G., Kombinatorische Anzahlbestimmungen f\"ur Gruppen,
Graphen und chemische Verbindungen, {\it Acta Math.} {\bf 68}
(1937), 145--254.

\bibitem{Pol}Polyakov A., Gauge f\/ields and space-time,
{\it Internat. J. Modern Phys.} {\bf 17} (2002), suppl., 119--136, \mbox{\href{http://arxiv.org/abs/hep-th/0110196}{hep-th/0110196}}.

\bibitem{Q1}Quillen D., Rational homotopy theory, {\it Ann. Math. (2)} {\bf 90}
(1969),  205--295.

\bibitem{Q2}Quillen D., Algebra cochains and cyclic cohomology, {\it Inst. Hautes \'Etudes Sci. Publ. Math.}
 (1988), no.~68, 139--174.

\bibitem{R} Ramis J.-P., S\'eries divergentes et th\'eories asymptotiques, {\it Bull. Soc. Math. France} {\bf  121} (1993),  suppl., 74~pages.

\bibitem{TZ}  Tradler T.,  Zeinallian M., On the cyclic Deligne conjecture,
{\it J. Pure Appl. Algebra} {\bf 204} (2006), 280--299,
\href{http://arxiv.org/abs/math.QA/0404218}{math.QA/0404218}.

\bibitem{Weil} Weil A., Elliptic functions according to Eisenstein and Kronecker,
{\it Ergebnisse der Mathematik und ihrer Grenzgebiete}, Vol.~88, Springer-Verlag, Berlin~-- New York, 1976.

\bibitem{WW} Whittaker E.T., Watson G.N., A course of modern analysis. An introduction to the general theory of inf\/inite processes and of analytic functions; with an account of the principal transcendental functions, reprint of 4th ed. (1927), Cambridge Mathematical Library, Cambridge University Press, Cambridge, 1996.

\bibitem{W}Witt E., Treue Darstellung Liescher Ringe, {\it J. Reine Angew. Math.} {\bf  177} (1937), 152--160.


\bibitem{Z} Zagier D., Elementary aspects of the Verlinde formula and of the Harder--Narasimhan--Atiyah--Bott formula, in  Proceedings of the Hirzebruch 65 Conference on Algebraic Geometry (Ramat Gan, 1993), {\it Israel Math. Conf. Proc.}, Vol.~9, Bar-Ilan Univ., Ramat Gan, 1996, 445--462.

\end{thebibliography}
\end{document}